# Orthogonal polynomial ensembles in probability theory*

**Wolfgang König**

*Universität Leipzig, Mathematisches Institut, Augustusplatz 10/11, D-04109, Leipzig, Germany*
*e-mail:* `koenig@math.uni-leipzig.de`

**Abstract:** We survey a number of models from physics, statistical mechanics, probability theory and combinatorics, which are each described in terms of an *orthogonal polynomial ensemble*. The most prominent example is apparently the Hermite ensemble, the eigenvalue distribution of the Gaussian Unitary Ensemble (GUE), and other well-known ensembles known in random matrix theory like the Laguerre ensemble for the spectrum of Wishart matrices. In recent years, a number of further interesting models were found to lead to orthogonal polynomial ensembles, among which the corner growth model, directed last passage percolation, the PNG droplet, non-colliding random processes, the length of the longest increasing subsequence of a random permutation, and others.

Much attention has been paid to universal classes of asymptotic behaviors of these models in the limit of large particle numbers, in particular the spacings between the particles and the fluctuation behavior of the largest particle. Computer simulations suggest that the connections go even farther and also comprise the zeros of the Riemann zeta function. The existing proofs require a substantial technical machinery and heavy tools from various parts of mathematics, in particular complex analysis, combinatorics and variational analysis. Particularly in the last decade, a number of fine results have been achieved, but it is obvious that a comprehensive and thorough understanding of the matter is still lacking. Hence, it seems an appropriate time to provide a surveying text on this research area.

In the present text, we introduce various models, explain the questions and problems, and point out the relations between the models. Furthermore, we concisely outline some elements of the proofs of some of the most important results. This text is aimed at non-experts with strong background in probability who want to achieve a quick survey over the field.

## 1. Introduction

In the 1950ies, it was found that certain important real $N$-particle ensembles (that is, joint distributions of $N$ real random objects) can be described by a

---

*This is an original survey paper

probability measure $\mathbb{P}_N$ of the form

$$\mathbb{P}_N(\mathrm{d}x) = \frac{1}{Z_N}\Delta_N(x)^2 \prod_{i=1}^{N} \mu(\mathrm{d}x_i), \qquad x = (x_1, \ldots, x_N) \in W_N, \qquad (1.1)$$

on the set

$$W_N = \{x = (x_1, \ldots, x_N) \in \mathbb{R}^N : x_1 < x_2 < \cdots < x_N\}, \qquad (1.2)$$

where $Z_N$ is the normalization, $\mu$ some distribution on $\mathbb{R}$ having all moments, and

$$\Delta_N(x) = \prod_{1 \leq i < j \leq N}(x_j - x_i) = \det\bigl[(x_i^{j-1})_{i,j=1,\ldots,N}\bigr], \qquad x = (x_1, \ldots, x_N),$$
(1.3)

is the well-known *Vandermonde determinant*. That is, $\mathbb{P}_N$ is the transformed configuration distribution of a vector of $N$ particles, distributed independently according to $\mu$ under the influence of the mutually repelling density $\Delta_N^2$, properly normalized to a probability measure on the so-called *Weyl chamber* $W_N$. The most important and one of the earliest examples is the joint distribution of the eigenvalues of a random matrix drawn from a Gaussian Unitary Ensemble (GUE), in which case $\mu$ is a Gaussian distribution, and $\mathbb{P}_N$ is called the *Hermite ensemble*. Also spectra of a couple of other types of random matrices turned out to admit a description of the form (1.1), among which the *Wishart matrices*, where $\mu$ is a Gamma distribution and $\mathbb{P}_N$ the *Laguerre ensemble*. The explicit form of (1.1) served as a starting point for many deep investigations of asymptotic spectral properties of random matrices. Furthermore, non-colliding Brownian motions (sometimes called Dyson's Brownian motions) could also successfully be investigated in the early 1960ies using descriptions in the spirit of (1.1). Also variants of (1.1) (e.g., with $\Delta_N^2$ replaced by $\Delta_N$ or by $\Delta_N^4$) turned out to have a significant relevance and could be treated using related methods.

For a long while, spectra distributions of certain random matrices (and the closely related non-colliding Brownian motions) were the only known important models that admit a description as in (1.1). However, in the second half of the 1990ies, the interest in non-colliding random processes was renewed and was put on a more systematic basis, and other types of statistical physics models were found to admit a description of the form (1.1): certain random growth models (equivalently, directed last passage percolation), polynuclear growth models, the problem of the length of the longest increasing subsequence in a random permutation, the Aztec diamond, and others. Furthermore, effective analytic techniques for deriving asymptotic properties of $\mathbb{P}_N$, which were developed in the early 1990ies, have recently been systematically extended and improved. As a consequence, in recent years a lot of questions about these models could be answered. The last ten years saw an exploding activity of research and an enormous progress in the rigorous understanding of some of the most important of these models, and the work is still going on with an increasing velocity. A significant number of deep and important results on universality questions have recently been solved, building on work of the last 40 or so years. However, it still seems as if a complete understanding of the phenomena has not yet been achieved, since many of the existing proofs are still based on explicit calculations and direct arguments. There seem some intricate mechanisms present which have been understood only in special cases by formal analogies. It will be an important and difficult task in future to find the essences of the phenomena in general.

In view of the great recent achievements, and also in order to draw the attention of non-experts to this field, it seems fruitful to write a comprehensive survey on most of the models that can be described by an ensemble as in (1.1). The present text is an attempt to explain the problems and questions of interest in a unifying manner, to present solutions that have been found, to give a flavor of the methods that have been used, and to provide useful guidelines to much of the relevant literature. It is aimed at the non-expert, the newcomer to the field, with a profound background in probability theory, who seeks a non-technical introduction, heuristic explanations, and a survey. Our concern is to comprehensively summarize the (in our opinion) most important available results and ideas, but not to lose ourselves in details or even technicalities. In the three remaining sections, we give an account on the three research areas we consider most important in connection with orthogonal polynomial ensembles: *random matrix theory*, *random growth models*, and *non-colliding random processes*.

A probability measure $\mathbb{P}_N$ of the form (1.1) is called an *orthogonal polynomial ensemble*. The theory of orthogonal polynomials is a classical subject, and appears in various parts of mathematics, like numerics, combinatorics, statistics and others. The standard reference on orthogonal polynomials is [Sz75]. However, the term 'orthogonal polynomial ensemble' is relatively recent and may be motivated by the following. Let $(\pi_N)_{N \in \mathbb{N}_0}$ denote the sequence of polynomials orthogonal with respect to the inner product on the space $L^2(\mu)$. The polynomials are unique by the requirement that the degree of $\pi_N$ is $N$, together with the normalization $\pi_N(x) = x^N + \mathcal{O}(x^{N-1})$. They may be obtained from the monomials $x \mapsto x^j$ via the well-known Gram-Schmidt algorithm. A nice relation[1] between the orthogonal polynomials and the ensemble $\mathbb{P}_N$ in (1.1) now is the fact that $\pi_N$ may be seen as the 'expected polynomial' of the form $\prod_{i=1}^{N}(x - x_i)$ with $(x_1, \ldots, x_N)$ distributed according to $\mathbb{P}_N$, i.e.,

$$\pi_N(x) = \int_{W_N} \prod_{i=1}^{N}(x - x_i)\, \mathbb{P}_N(\mathrm{d}x_1 \cdots \mathrm{d}x_N). \tag{1.4}$$

## 2. Random matrix theory

In spite of the appearance of various random matrix distributions in several areas of mathematics and physics, it has become common to use the term *random matrix theory* exclusively for those matrix distributions that are used, since

---

[1] Further connections will be exploited in Section 2.7 below.

Wigner's introduction to physics in the 1950ies, as models for energy levels in slow nuclear reactions. Measurements had already given rise to a hope that the energy levels follow a universal picture. Wigner's hope was that the eigenvalues of appropriate classes of random matrices would be mathematically tractable and would reflect this universality in the limit of unbounded matrix size. Based on Wigner's early work, Dyson [Dy62a], [Dy62c] argued on physical grounds that three certain matrix classes be relevant for the description of energy levels, the by now famous orthogonal, unitary and symplectic *Gaussian ensembles*. It soon turned out that their eigenvalue distributions are given in terms of certain orthogonal polynomial ensembles. In the mid-nineties, seven more Gaussian random matrix ensembles were introduced [Ve94], [AZ96], [AZ97], [Zi97], and it was argued that these in total ten classes form a complete classification of the set of random matrix ensembles that are physically relevant in a certain sense, subject to some symmetry constraints.

In the last decades, random matrix theory became a major mathematical and physical research topic, and more and more exciting phenomena were found. In the last decade, universality of many aspects could be proven for large classes of models, and the research is going on increasingly fast.

The standard reference on the mathematical treatment of random matrices is [Me91]. Authored by a physicist with strong mathematical interest, it explains the physical relevance of a host of random matrix models and provides a great amount of relevant formulas and calculations. A recent historical survey on the field from a physicist's point of view in [FSV03] (see the entire volume), which contains a vast list of references, mostly from the physics literature. A thorough summary of the proofs of some of the most important results on random matrix theory from the viewpoint of Riemann-Hilbert theory is in [De98]. Further surveying and appetizing texts on random matrix theory are [TW93b] and [Di03]. When the present text is being written, some (teams of) authors are preparing monographs on random matrix theory, among which [Fo05+].

In the present section we first introduce to some of the above mentioned matrix ensembles and their eigenvalue distributions in Sections 2.1–2.4, present the famous Wigner semicircle law in Section 2.5, discuss correlation functions in Section 2.6 and introduce the important method of orthogonal polynomials in Section 2.7. Afterwards, we present the most important asymptotic results on eigenvalues, the bulk limits in Section 2.8 and the edge asymptotics in Section 2.9. The main proof method, the Riemann-Hilbert theory, is outlined in Section 2.10. Finally, in Section 2.11 we explain some relations to the zeros of the famous Riemann zeta function.

### 2.1 The questions under interest

Consider a random Hermitian $N \times N$-matrix, $M$, and denote its eigenvalues by $\lambda_1 \leq \lambda_2 \leq \cdots \leq \lambda_N$. Hence, $\lambda = (\lambda_1, \ldots, \lambda_N)$ is a random element of the closure of the Weyl chamber $W_N$ in (1.2). Among others, we shall ask the following questions:

(i) What is the distribution of $\lambda$ for fixed $N \in \mathbb{N}$?

  (ii) What is the limiting scaled distribution of $\lambda$ as $N \to \infty$, in terms of the empirical measure $\frac{1}{N}\sum_{i=1}^{N} \delta_{\widetilde{\lambda}_i}$, for an appropriate scaling $\widetilde{\lambda}_i$ of $\lambda_i$?
  (iii) What is the limiting behavior of the largest eigenvalue, $\lambda_N$, as $N \to \infty$? (Or of the smallest, $\lambda_1$, or the joint distribution of a few of the smallest, say $(\lambda_1, \ldots, \lambda_m)$ for some $m$.) More precisely, what is the right normalization for a law of large numbers, and what is the right scaling for a limit law, if present?
  (iv) What are the limiting statistics of the spacings between neighboring eigenvalues? How many gaps are there with a given maximal length? What is the average distance between $\lambda_{cN-r_N}$ and $\lambda_{cN+r_N}$ for some $c \in (0,1)$ and some $r_N \to \infty$ such that $r_N/N \to \infty$?

Question (iii) refers to the *edge of the spectrum*, while (iv) refers to the *bulk of the spectrum*.

The so-called *Wigner surmise* conjectures that the limiting spacing between two subsequent eigenvalues of a large Gaussian matrix should have the density $(0,\infty) \ni x \mapsto Cxe^{-cx^2}$. This is true for a $(2\times 2)$-matrix $\left(\begin{smallmatrix} a & b \\ b & c \end{smallmatrix}\right)$ with independent standard Gaussian entries $a, b, c$: The spacing $\lambda_2 - \lambda_1$ is equal to $[(a-c)^2 + 4b^2]^{1/2}$, whose square has the $\chi^2$-distribution. However, the Wigner surmise turned out to be inaccurate (even though rather close to the true distribution): the asymptotic spacing distribution is different.

**2.2 Matrix distributions**

It turned out [Dy62a] that, according to time reversal invariance properties of the material considered, basically three different matrix classes are of interest as models for energy levels of nuclea: matrices whose entries are (1) real numbers, (2) complex numbers, and (3) quaternions. One basic requirement is that the random matrices considered be symmetric, respectively Hermitian, respectively self-dual, such that all the eigenvalues are real numbers. For the (famous and most studied) special case of Gaussian entries, these three cases correspond to the *Gaussian Orthogonal Ensemble* (GOE), the *Gaussian Unitary Ensemble* (GUE) and the *Gaussian Simplectic Ensemble* (GSE). In the following, we shall concentrate mostly on the unitary ensemble, since this class is, in some respects, technically the easiest to treat and exhibits the farthest reaching connections to other models.

We assume that $M = (M_{i,j})_{i,j=1,\ldots,N}$ is a random Hermitian $(N \times N)$-matrix with complex entries. In particular, the diagonal entries $M_{i,i}$ are real, and for $i \neq j$, we have $M_{i,j} = M_{i,j}^{(\mathrm{R})} + \mathrm{i}\, M_{i,j}^{(\mathrm{I})} = M_{j,i}^{(\mathrm{R})} - \mathrm{i}\, M_{j,i}^{(\mathrm{I})} = \overline{M}_{j,i}$, where $M_{j,i}^{(\mathrm{R})}$ and $M_{j,i}^{(\mathrm{I})}$ are the real part and imaginary part, respectively.

Two basic respective requirements are (1) independence of the matrix entries, and (2) invariance of the distribution of the matrix under unitary conjugations. These two ideas lead to different matrix classes:

**Wigner matrices:** *We call the random Hermitian matrix $M$ a Wigner matrix if the collection $\{M_{i,j}^{(\mathrm{R})}: i,j = 1,\ldots,N,\ i < j\} \cup \{M_{i,j}^{(\mathrm{I})}: i,j = 1,\ldots,N,\ i < j\} \cup \{M_{i,i}: i = 1,\ldots,N\}$ consists of independent, not necessarily identically* *distributed, random variables with mean zero and a fixed positive variance, which is the same for the real parts and for the imaginary parts, respectively.*

Hence, there are $N^2$ independent real variables that determine the distribution of $M$. The distribution of the diagonal elements is arbitrary, subject to moment conditions.

**Unitary-invariant matrices:** *We call the random Hermitian matrix $M$ unitary-invariant if the joint distribution of its entries has the form*

$$\mathbb{P}(\mathrm{d}M) = \text{const. } e^{-F(M)} \prod_{i=1}^{N} \mathrm{d}M_{i,i} \prod_{1 \leq i < j \leq N} \left[ \mathrm{d}M_{i,j}^{(\mathrm{R})} \, \mathrm{d}M_{i,j}^{(\mathrm{I})} \right] = \text{const. } e^{-F(M)} \, \mathrm{d}M, \tag{2.1}$$

*for some function $F$, and, for every unitary matrix $U$, $\mathbb{P}(\mathrm{d}M)$ is invariant under conjugation with $U$.*

The last requirement implies that $e^{-F(UMU^{-1})} \, \mathrm{d}[UMU^{-1}] = e^{-F(M)} \, \mathrm{d}M$, and since it is easy to see that $\mathrm{d}[UMU^{-1}] = \mathrm{d}M$ [De98, p. 92], the function $F$ must be a symmetric function of the eigenvalues of $M$. One particular and important example is the case where

$$F(M) = \mathrm{Tr}(Q(M)), \qquad Q(x) = \gamma_{2j} x^{2j} + \cdots + \gamma_0, \text{a polynomial with } \gamma_{2j} > 0. \tag{2.2}$$

With the exception of the Gaussian case $j = 1$, there are strong correlations between all the matrix entries. The idea behind the invariance under unitary conjugations is that the matrix distribution should not depend on the observation system, as long as it is based on a unitary coordinate axis.

The famous GUE lies in the intersection of the Wigner-class and the unitary-invariant class. It is a Wigner matrix with all the sub-diagonal entries being complex standard normal variables[2] and the diagonal entries being a real normal variable with variance two[3]. Alternately, it is the unitary-invariant matrix of the form (2.1) with $F(M) = \mathrm{Tr}(M^2)$.

The GOE is the real variant of the GUE; i.e., the sub-diagonal entries are independent standard real normal variables with the same variance as the diagonal entries. Hence, the GOE has $\frac{1}{2}N(N+1)$ independent sources of real randomness.

The GSE is the symplectic variant of the GUE, i.e., the diagonal entries are real standard normals as in the GUE, and the sub-diagonal entries are elements of the quaternion numbers. Their four components are i.i.d. real standard normal variables. Hence, the GSE has $N + 2N(N-1)$ independent real randomnesses.

Further important related classes of random matrices are the *Wishart matrices*, which are of the form $A^*A$ with $A$ a (not necessarily square) matrix having throughout i.i.d. complex normal entries (first considered in a multivariate statistics context by Wishart [Wi28]). See [Me91] for further classes.

---

[2]By this we mean that the real and the imaginary part are two independent standard normal variables.

[3]Some authors require the sum of the variances per entry to be equal to one, or equal to $1/N$.

### 2.3 Eigenvalue distributions

Let $\lambda_1 \leq \lambda_2 \leq \cdots \leq \lambda_N$ be the $N$ eigenvalues of the random Hermitian matrix $M$. We ask for the distribution of the random vector $\lambda = (\lambda_1, \ldots, \lambda_N)$. A concise answer for a general Wigner matrix $M$ seems inaccessible, but for unitary-invariant ensembles there is a nice, fundamental formula. We formulate the GUE case and make a couple of remarks afterwards.

**Lemma 2.1 (Eigenvalue distribution for GUE).** *Let $M$ be a random matrix from GUE. Then the distribution of the vector $\lambda = (\lambda_1, \ldots, \lambda_N)$ of eigenvalues of $M$ has the density*

$$P_N(x) = \frac{1}{Z_N} \Delta_N(x)^2 \prod_{i=1}^N e^{-x_i^2}, \qquad x = (x_1, \ldots, x_N), \tag{2.3}$$

*with $Z_N$ the appropriate normalizing constant on the Weyl chamber $W_N$ in (1.2).*

**Sketch of the proof.** Choose a (random) unitary matrix $U$ which diagonalizes $M$, i.e., the matrix $D = UMU^{-1}$ is the diagonal matrix with the eigenvalues on the main diagonal. Hence,

$$\begin{aligned} \mathrm{d}M &= \mathrm{d}(U^* D U) = \mathrm{d}U^* \cdot D \cdot U + U^* \cdot \mathrm{d}D \cdot U + U^* \cdot D \cdot \mathrm{d}U \\ &= U^* \cdot (\mathrm{d}D + U \cdot \mathrm{d}U^* \cdot D + D \cdot \mathrm{d}U \cdot U^*) \cdot U \\ &= \mathrm{d}D + U \cdot \mathrm{d}U^* \cdot D + D \cdot \mathrm{d}U \cdot U^* \\ &= \mathrm{d}D + \mathrm{d}A \cdot D - D \cdot \mathrm{d}A, \end{aligned} \tag{2.4}$$

where we used the invariance of $\mathrm{d}M$ under unitary conjugations, and we introduced $\mathrm{d}A = U \cdot \mathrm{d}U^* = -\mathrm{d}U \cdot U^*$. Now integrate over $\mathrm{d}M_{i,j}$ with $i < j$ and use calculus. See [Me91, Ch. 3] or [HP00, Ch. 4] for details. □

**Remark 2.2.** (i) We chose the normalization $Z_N$ such that $P_N$ is normalized on $W_N = \{x \in \mathbb{R}^N : x_1 < x_2 < \cdots < x_N\}$. We extend $P_N$ to a permutation symmetric function on $\mathbb{R}^N$. Hence, $\widehat{P}_N = \frac{1}{N!} P_N$ is a probability density on $\mathbb{R}^N$.

(ii) The density in (2.3) is called the *Hermite ensemble*. This is one of the most prominent examples of an orthogonal polynomial ensemble; the name refers to the Hermite polynomials which form an orthonormal base with respect to Gaussian weights.

(iii) For the GOE and the GSE, there are analogous formulas. Indeed, replace $\Delta_N(x)^2$ by $\Delta_N(x)^\beta$ with $\beta = 1$ respectively $\beta = 4$ to obtain the corresponding statement for GOE, respectively for GSE [Me91, Ch. 3]. The three matrix classes are usually marked by the parameter

$$\beta = 1 \text{ for GOE}, \qquad \beta = 2 \text{ for GUE}, \qquad \beta = 4 \text{ for GSE}. \tag{2.5}$$

(iv) It is easy to extend Lemma 2.1 to unitary-invariant matrix distributions. Indeed, if the distribution of $M$ is of the form (2.1) with $F(M) = f(\lambda_1, \ldots, \lambda_N)$, a symmetric function that depends only on the spectrum of $M$, then the density of $(\lambda_1, \ldots, \lambda_N)$ is proportional to $x \mapsto e^{-f(x)} \Delta_N(x)^2$. An analogous assertion is true for the orthogonal case, see [HP00, Ch. 4].

(v) If $M = A^* A$ is a Wishart matrix, i.e., $A$ is an $(N \times k)$-matrix (with $k \leq N$) having throughout independent complex standard normal entries, then the vector of eigenvalues of $M$ has the density [Ja64]

$$x \mapsto \frac{1}{Z_{N,k}} \Delta_N(x)^2 \prod_{i=1}^N [x_i^{N-k} e^{-x_i}], \qquad x \in W_N \cap (0, \infty)^N. \quad (2.6)$$

This ensemble is called the *Laguerre* ensemble.

(vi) Using *Selberg's integral* [HP00, p. 118/9], the normalizing constants of the Hermite ensemble and the Laguerre ensemble may be identified in terms of the Gamma-function. Indeed, for any $\beta > 0$, we have

$$\int_{\mathbb{R}^N} |\Delta_N(x)|^\beta \frac{e^{-\|x\|^2}}{(2\pi)^{N/2}} \, dx = \prod_{i=1}^N \frac{\Gamma(1 + i\frac{\beta}{2})}{\Gamma(1 + \frac{\beta}{2})}, \quad (2.7)$$

and, for any $a > 0$,

$$\int_{\mathbb{R}^N} |\Delta_N(x)|^\beta \prod_{i=1}^N [x_i^{a-1} e^{-x_i}] \, dx = \prod_{j=0}^{N-1} \frac{\Gamma(1 + (1+j)\frac{\beta}{2}) \Gamma(a + j\frac{\beta}{2})}{\Gamma(1 + \frac{\beta}{2})}. \quad (2.8)$$

(vii) There is obviously a mutually repelling force between the eigenvalues in (2.3): the density vanishes if any two of the $N$ arguments approach each other. It does not seem easy to derive an intuitive reason for this repellence from random matrix considerations, but if the matrix $M$ is embedded in a natural process of random Hermitian matrices, then the process of eigenvalues admits a nice identification that makes the repellence rather natural. This is the subject of Section 4.1 below.  ◇

### 2.4 Circular ensembles

An important type of random Gaussian matrices are the *circular ensembles*, which were introduced in [Dy62a] in the desire to define a matrix model that can be seen as the conditional Gaussian ensembles given a fixed value of the exponential weight $F(M)$ in (2.1). Again, there is an orthogonal, unitary and symplectic version of the circular ensemble.

We give the definition of the circular ensembles [Me91, Ch. 9]. The *circular orthogonal ensemble (COE)* is the unique distribution on the set of orthogonal symmetric $(N \times N)$-matrices that is invariant under conjugation with any real orthogonal matrix. That is, an orthogonal symmetric random matrix $S$ is COE-distributed if and only if $WSW^{-1}$ has the same distribution as $S$, for any real orthogonal matrix $W$. The *circular unitary ensemble (CUE)* is the unique distribution on the set of complex unitary $(N \times N)$-matrices that is invariant under (two-sided) transformations with unitary matrices, i.e., a complex unitary random matrix $S$ is CUE-distributed if and only if $USV$ has the same distribution as $S$, for any two unitary matrices $U$ and $V$. Finally, the *circular symplectic ensemble (CSE)* is the unique distribution on the set of self-dual unitary quaternion matrices that is invariant under every automorphism $S \mapsto W^{\text{R}}SW$, where $W$ is any unitary quaternion matrix and $W^{\text{R}}$ its dual.

All eigenvalues of the circular matrices lie on the unit circle and may be written $\lambda_i = e^{\mathrm{i}\,\theta_i}$ with $0 \leq \theta_1 < \theta_2 < \cdots < \theta_N < 2\pi$. One advantage of the circular ensembles is that the joint distribution density of their eigenvalues admits particularly simple formulas. Indeed, adopting the parameter $\beta = 1, 2, 4$ for the COE, CUE and CSE, respectively (recall (2.5)), the density of the vector $(\theta_1, \ldots, \theta_N)$ of eigenvalue angles is given as

$$P_N^{(\text{circ},\beta)}(\theta_1, \ldots, \theta_N) = \frac{1}{Z_N^{(\text{circ},\beta)}} \prod_{1 \leq \ell < j \leq N} |e^{\mathrm{i}\,\theta_\ell} - e^{\mathrm{i}\,\theta_j}|^\beta = \frac{1}{Z_N^{(\text{circ},\beta)}} |\Delta_N(e^{\mathrm{i}\,\theta_\cdot})|^\beta. \tag{2.9}$$

Here we chose the normalization such that $P_N^{(\text{circ},\beta)}$ is a probability density on $W_N \cap [0, 2\pi)^N$ where $W_N = \{x \in \mathbb{R}^N \colon x_1 < \cdots < x_N\}$ is the Weyl chamber.

### 2.5 The law of large numbers: Wigner's semi-circle law

In this section we present the famous semi-circle law first proved by Wigner: the convergence of the mean eigenvalue density as the size of the matrix increases to infinity. This is an asymptotic statement about the convergence of the empirical measure[4] of the appropriately scaled eigenvalues of a random matrix towards the distribution

$$\frac{\mu_*(\mathrm{d}x)}{\mathrm{d}x} = \frac{1}{\pi}\sqrt{2 - x^2}\,\mathbb{1}_{[-\sqrt{2},\sqrt{2}]}(x), \tag{2.10}$$

the famous *semicircle distribution*. We first formulate the semicircle law for the GUE, make some remarks and sketch two proofs. Afterwards we summarize some extensions.

**Theorem 2.3 (Semicircle law).** *Let the random matrix $M_N$ be a GUE-matrix of size $N$, with eigenvalues $\lambda_1^{(N)} < \cdots < \lambda_N^{(N)}$. Let*

$$\mu_N = \frac{1}{N}\sum_{i=1}^N \delta_{\widetilde{\lambda}_i^{(N)}}, \qquad \text{where } \widetilde{\lambda}_i^{(N)} = N^{-\frac{1}{2}}\lambda_i^{(N)}, \tag{2.11}$$

*be the empirical measure of the rescaled eigenvalues. Then $\mu_N$ converges weakly in distribution towards the semicircle distribution $\mu_*$ in (2.10).*

We shall call $\lambda_1^{(N)}, \ldots, \lambda_N^{(N)}$ the (unscaled) eigenvalues and $\widetilde{\lambda}_1^{(N)}, \ldots, \widetilde{\lambda}_N^{(N)}$ the (re)scaled eigenvalues.

---

[4]By the *empirical measure* of $N$ points $x_1, \ldots, x_N$ we denote the probability measure $\frac{1}{N}\sum_{i=1}^N \delta_{x_i}$.

**Remark 2.4.** (i) Theorem 2.3 reveals that the eigenvalues are of order $\sqrt{N}$ and that the largest behaves like $\sqrt{2N}$. In particular, all eigenvalues lie in the interval $N^{1/2}[-\sqrt{2}-\varepsilon, \sqrt{2}+\varepsilon]$ for any $\varepsilon > 0$ with overwhelming probability, and the spacings between subsequent eigenvalues are of order $N^{-1/2}$ in the bulk of the spectrum and much larger close to the edge.
 (ii) The convergence takes place in the sense that the expectation of every bounded and continuous function of the empirical measure converges. Note that the moments (i.e., the family of maps $\mu \mapsto \int x^k \mu(\mathrm{d}x)$ for $k \in \mathbb{N}$) constitute a convergence determining family.
(iii) Note that, for any $a < b$,

$$\mathbb{E}\big[\mu_N([a,b])\big] = \frac{1}{N}\mathbb{E}\big[\#\{i\colon \widetilde{\lambda}_i \in [a,b]\}\big].$$

In particular, the semicircle law states that the expected number of unscaled eigenvalues $\lambda_i^{(N)}$ in the interval $[aN^{\frac{1}{2}}, bN^{\frac{1}{2}}]$ behaves like $N\mu_*([a,b])$. See Remark 2.5(vi) for further asymptotic statements.
(iv) The convergence in Theorem 2.3 has also been proved [Ar67] in distribution in the almost sure sense, see [HP00, Th. 4.1.5]. More precisely, let $(M_{i,j})_{i,j\in\mathbb{N}}$ be a sequence of independent standard complex normal random variables and denote by $M_N = (M_{i,j})_{i,j\leq N}$ the $(N\times N)$-corner. Let $\mu_N$ (as in (2.11)) denote the empirical measure of the rescaled eigenvalues of $M_N$. Then all the $k$-th moments of $\mu_N$ converge towards the $k$-th moment of $\mu_*$.
 (v) See [HP00, Ch. 4] for the statement analogous to Theorem 2.3 for the orthogonal ensembles. ◇

We turn now to sketchs of two proofs.

**Sketch of the first proof: the method of moments.** This is Wigner's original method [Wi55], [Wi58], see [HP00, Ch. 4]. The idea is that it suffices to prove that the expected moments of $\mu_N$ converge to the ones of $\mu_*$, i.e.,

$$\lim_{N\to\infty} \mathbb{E}\Big[\int_\mathbb{R} x^k \mu_N(\mathrm{d}x)\Big] = \int_\mathbb{R} x^k \mu_*(\mathrm{d}x), \qquad k \in \mathbb{N}. \tag{2.12}$$

By symmetry, all odd moments of both $\mu_N$ and $\mu_*$ are zero, hence it suffices to consider $k = 2m$. The $(2m)$-th moments of $\mu_*$ are known to be $\frac{2^{-m}}{1+m}\binom{2m}{m}$. Note that the left hand side is equal to the normalized trace of $M_N^{2m}$, i.e.,

$$\mathbb{E}\Big[\int_\mathbb{R} x^{2m} \mu_N(\mathrm{d}x)\Big] = \frac{1}{N}\sum_{i=1}^N \mathbb{E}\Big[\int x^{2m} \delta_{\widetilde{\lambda}_i^{(N)}}(\mathrm{d}x)\Big] = \frac{1}{N^{1+m}}\mathbb{E}\Big[\sum_{i=1}^N \big(\lambda_i^{(N)}\big)^{2m}\Big]$$

$$= \frac{1}{N^{1+m}}\mathbb{E}\big[\mathrm{Tr}(M_N^{2m})\big] = \frac{1}{N^{1+m}} \sum_{i_1,\ldots,i_{2m}=1}^N \mathbb{E}\Big[\prod_{j=1}^{2m} M_{i_{j-1}, i_j}\Big],$$
(2.13)

where $M_{i,k}$ denote the entries of the matrix $M_N$. Some combinatorial work has to be done in order to discard from the sum those terms that do not contribute, and to extract the leading terms, using the independence of the matrix entries and rough bounds on the moments of the matrix entries. The term coming from the subsum over those multi-indices $i_1, \ldots, i_{2m}$ with $\#\{i_1, \ldots, i_{2m}\} < m+1$ is shown to vanish asymptotically, and the one with $\#\{i_1, \ldots, i_{2m}\} > m+1$ is shown to be equal to zero. □

The second proof is in the spirit of statistical mechanics and is based on the eigenvalue density in (2.3). Indeed, the convergence is derived with the help of large-deviation type arguments and the minimization of a certain energy functional. In particular, the semicircle law turns out to be the unique minimizer, because of which it is sometimes called an *equilibrium measure* for that functional. We partially follow the presentation in [De98, Ch. 6], which is based on [Jo98] and [DMK98]. A general reference for equilibrium measures and related material is [ST97].

**Sketch of the second proof: the equilibrium measure method.** The starting point is the observation that the joint density $P_N$ of the unscaled eigenvalues in (2.3) is of the form $P_N(x) = \frac{1}{Z_N} e^{-H_N(x)}$ with the Hamiltonian

$$H_N(x) = \sum_{i=1}^{N} x_i^2 - 2 \sum_{1 \leq i < j \leq N} \log(x_j - x_i). \tag{2.14}$$

In order to obtain a non-degenerate limit law, we have to rescale the $\lambda_i^{(N)}$ in such a way that both parts of $H_N(x)$ are of the same order in $N$. Since the second part is always of order $N^2$, it is clear that we should consider the scaling $\widetilde{\lambda}_i^{(N)} = N^{-\frac{1}{2}} \lambda_i^{(N)}$ as in the theorem. The vector $\widetilde{\lambda}^{(N)}$ of the rescaled quantities has the density

$$\mathbb{P}(\widetilde{\lambda}^{(N)} \in dx) = \frac{1}{\widetilde{Z}_N} e^{-N^2 \widetilde{H}_N(x)} \, dx, \tag{2.15}$$

where

$$\widetilde{H}_N(x) = \frac{1}{N} \sum_{i=1}^{N} x_i^2 - \frac{2}{N^2} \sum_{1 \leq i < j \leq N} \log(x_j - x_i), \tag{2.16}$$

and we absorbed some terms in the new normalizing constant. In terms of the empirical measure of the rescaled quantities, $\mu_N$, the Hamiltonian takes the shape $\widetilde{H}_N \approx \mathcal{I}(\mu_N)$, where

$$\mathcal{I}(\mu) = \int_{\mathbb{R}} x^2 \, \mu(dx) - \int_{\mathbb{R}} \int_{\mathbb{R}} \log|x - y| \, \mu(dx) \mu(dy), \qquad \mu \in \mathcal{M}_1(\mathbb{R}). \tag{2.17}$$

Here we suppressed the diagonal terms, i.e., the summands for $i = j$, which is a technical issue. Since the integration is only of the order $N$ and the exponent of order $N^2$, it is clear that the large-$N$ behavior of the measure $\frac{1}{\widetilde{Z}_N} e^{-N^2 \mathcal{I}(\mu_N)} \, d^N x$ is determined by the minimizer(s) of the variational problem

$$E = \inf_{\mu \in \mathcal{M}_1(\mathbb{R})} \mathcal{I}(\mu). \tag{2.18}$$

The minimizer(s) are called *equilibrium measure(s)*. According to (a high-dimensional variant of) the well-known Laplace method, the value of $E$ should be the large-$N$ exponential rate of $P_N(x)$ with speed $N^2$, and the empirical measures $\mu_N$ should converge towards the minimizer(s).

The analysis of the problem in (2.18) is not particularly difficult. Using standard methods one shows the existence and uniqueness of the equilibrium measure and the compactness of its support. Using the Euler-Lagrange equation in the interior of its support, one identifies the equilibrium measure with the semicircle law, $\mu_*$. However, in order to show the convergence of $\mu_N$ towards $\mu_*$, one needs to show that the contribution coming from outside a neighborhood of $\mu_*$ is negligible, which is a more difficult issue. This is carried out in [Jo98]. □

**Remark 2.5.** (i) The moment method has been generalized to a large class of Wigner matrices, [HP00, Ch. 4]. Indeed, assume that the real and the imaginary parts of the sub-diagonal entries of $M_N$ are independent, centred and have variance one and that the diagonal entries have variance two, and assume that, for any $k \in \mathbb{N}$, their $k$-th moments are uniformly bounded in $N$, then the conclusion of Theorem 2.3 remains true.

(ii) The equilibrium measure method has been generalized to a large class of unitary-invariant matrices in [D99] and [DMK98], see also [De98]. To mention one of the most obvious generalisations, let $M_N$ be a matrix as in (2.1) with $F$ as in (2.2), i.e., the eigenvalues have the density in (2.3) with the term $x_i^2$ replaced by the polynomial $Q(x_i) = \gamma_{2j} x_i^{2j} + \mathcal{O}(x_i^{2j-1})$; recall Remark 2.2(iv). The correct scaling is $\widetilde{\lambda}_i^{(N)} = N^{-\frac{1}{2j}} \lambda_i^{(N)}$, and in the limit as $N \to \infty$, only the leading term of $Q(x_i)$ survives. The empirical measure of the $\widetilde{\lambda}_i^{(N)}$ converges weakly towards the equilibrium measure of the functional

$$\mu \mapsto \int_{\mathbb{R}} \gamma_{2j} x^{2j} \, \mu(\mathrm{d}x) - \int_{\mathbb{R}} \int_{\mathbb{R}} \log|x-y| \, \mu(\mathrm{d}x) \mu(\mathrm{d}y). \qquad (2.19)$$

The analysis of this functional and the proof of convergence towards its minimizer is similar to the proof in the special case where $Q(x) = x^2$. The equilibrium measure has a density, and its support is compact. If $\psi$ denotes the density and $[-a, a]$ its support, then $\psi(x) = (a^2 - x^2)^{\frac{1}{2}} h_1(x)$ for $|x| < a$, where $h_1$ is a polynomial of order $2j - 2$.

(iii) Even more generally, one starts immediately from distributions as in (2.3) with $x_i^2$ replaced by $NV(x_i)$ (note the factor of $N$) with some sufficiently regular function $V$ tending to infinity at infinity sufficiently fast. With this ansatz, no rescaling is necessary, i.e., the empirical measure of the unscaled vector $(\lambda_1^{(N)}, \ldots, \lambda_N^{(N)})$ converges. The relevant functional is then the one in (2.19) with $\gamma_{2j} x^{2j}$ replaced by $V(x)$. The Euler-Lagrange equations for this functional are, for some $l \in \mathbb{R}$,

$$2\int_{\mathbb{R}} \log|x-y|^{-1} \mu_*(\mathrm{d}y) + V(x) = l \qquad x \in \mathrm{supp}(\mu_*)^\circ. \qquad (2.20)$$

However, for general $V$, the explicit identification of the minimizer is considerably more difficult and involved. In general, if $V$ is convex, then the support of the equilibrium measure is still an interval, but in the general case it consists of a finite union of intervals, provided that $V$ is analytic [DMK98].

(iv) The energy functional $\mathcal{I}$ in (2.17) has an interpretation in terms of an electrostatic repulsion in the presence of an external quadratic field, if $\mu$ is the distribution of electrons. The second term is sometimes called *logarithmic entropy* or *Voiculescu's entropy*, see [Vo93] and [Vo94].

(v) An advantage of the equilibrium measure method is that it opens up the possibility of a large-deviation principle for the empirical measure of the rescaled eigenvalues. (This is, roughly speaking, the determination of the large-$N$ decay rate of the probability for a deviation of the empirical measure from the semicircle law in terms of a variational problem involving the energy functional.) The first proof of such a principle is in [BAG97], after pioneering (and less rigorous) work in [Vo93] and [Vo94]. Extensive and accessible lecture notes on large deviation techniques for large random matrices may be found in [Gui04].

(vi) In the course of the equilibrium-measure proof of Theorem 2.3 (see [De98, Theorem 6.96]), for every $k \in \mathbb{N}$, also the weak convergence of the $k$-dimensional marginal of $\widehat{P}_N$ with density

$$\widehat{P}_{N,k}(x_1,\ldots,x_k) = \Big(\int_{\mathbb{R}^{N-k}} \widehat{P}_N(x_1,\ldots,x_N)\,\mathrm{d}x_{k+1}\cdots\mathrm{d}x_N\Big), \qquad (2.21)$$

towards the $k$-fold product measure $\mu_*^{\otimes k}$ is proved. As an elementary consequence, $N^{-k}$ times the expected number of $k$-vectors of different rescaled eigenvalues in $[a,b]$ converges towards $\mu_*([a,b])^k$.

(vii) There is an analogue of the semicircle law for the spectrum of the circular ensembles introduced in Section 2.4, without normalisation of the eigenvalues required. An innovative technique for deriving this law was introduced in [DS94] (see also [DE01]), where the asymptotic independence and normality of the traces of powers of the random matrix under consideration is shown. Related results are derived in [DS94] for the problem of the longest increasing subsequence of a uniform random permutation, which is introduced in Section 3.5.

$\diamond$

## 2.6 Correlation functions

In this section we let $P_N\colon W_N \to [0,\infty)$ be any probability density on the Weyl chamber $W_N$ in (1.2) and $\lambda = (\lambda_1,\ldots,\lambda_N) \in W_N$ be a random variable with density $P_N$. We introduce the so-called correlation functions of $P_N$, which will turn out to be important for two reasons: (1) much interesting information about the random variable $\lambda$ can be expressed in terms of the correlation functions, and (2) when specializing $P_N$ to an orthogonal polynomial ensemble, the correlation functions admit a determinantal representation which will be fundamental for the asymptotic analysis of the ensemble.

We extend $P_N$ to a permutation invariant function on $\mathbb{R}^N$. Then $\widehat{P}_N = \frac{1}{N!} P_N$ is a probability density on $\mathbb{R}^N$. For $k \in \mathbb{N}$, the *k-point correlation function* is defined as

$$R_k^{(N)}(x_1, \ldots, x_k) = \frac{N!}{(N-k)!} \int_{\mathbb{R}^{N-k}} \widehat{P}_N(x) \, \mathrm{d}x_{k+1}\mathrm{d}x_{k+2}\cdots \mathrm{d}x_N, \qquad x_1, \ldots, x_k \in \mathbb{R}. \tag{2.22}$$

Then $\widehat{P}_{N,k} = \frac{(N-k)!}{N!} R_k^{(N)}$ is a probability density on $\mathbb{R}^k$, the marginal density of $\widehat{P}_N$ in (2.21). It is a simple combinatorial exercise to see that, for any measurable set $A \subset \mathbb{R}$, the quantity $\int_{A^k} R_k^{(N)}(x) \, \mathrm{d}^k x$ is equal to the expected number of $k$-tuples $(\lambda_{i_1}, \ldots, \lambda_{i_k})$ of distinct particles such that $\lambda_{i_j} \in A$ for all $j = 1, \ldots, k$. In particular, $R_1^{(N)}(x) \, \mathrm{d}x$ is the expected number of particles in $\mathrm{d}x$.

As a first important application, the probability that a given number of particles lie in a given set can be expressed in terms of the correlation functions as follows.

**Lemma 2.6.** *For any $N \in \mathbb{N}$, any $m \in \{0, 1, \ldots, N\}$ and any interval $I \subset \mathbb{R}$,*

$$\mathbb{P}\big(\#\{i \leq N : \lambda_i \in I\} = m\big) = \frac{(-1)^m}{m!} \sum_{k=m}^{N} \frac{(-1)^k}{(k-m)!} \int_{I^k} R_k^{(N)}(x) \, \mathrm{d}^k x. \tag{2.23}$$

**Sketch of the proof.** We only treat the case $m = 0$, the general case being a simple extension. Expand

$$\prod_{i=1}^{N} [1 - \mathbb{1}_I(x_i)] = \sum_{k=0}^{N} (-1)^k \zeta_k\big(\mathbb{1}_I(x_1), \ldots, \mathbb{1}_I(x_N)\big)$$

where the functions $\zeta_k$ are permutation symmetric polynomials, which are defined by the relation $\prod_{i=1}^{N}(z - \alpha_i) = \sum_{k=0}^{N}(-1)^{N-k} z^k \zeta_k(\alpha)$ for any $z \in \mathbb{R}$ and $\alpha = (\alpha_1, \ldots, \alpha_N) \in \mathbb{R}^N$. Now multiply by the density $\widehat{P}_N$ and integrate over $\mathbb{R}^N$. Using the explicitly known coefficients of the polynomials $\zeta_k$, and using the permutation invariance of $R_k^{(N)}$, one arrives at (2.23). □

Also the expected number of spacings in the vector $\lambda$ can be expressed in terms of the correlation functions. For $x = (x_1, \ldots, x_N) \in W_N$, $u \in \mathbb{R}$ and $s, t \geq 0$ denote by

$$S^{(N)}(s; x) = \#\big\{j \in \{1, \ldots, N-1\} : x_{j+1} - x_j \leq s\big\}, \tag{2.24}$$

$$S_t^{(N)}(s, u; x) = \#\big\{j \in \{1, \ldots, N-1\} : x_{j+1} - x_j \leq s, |x_j - u| \leq t\big\}, \tag{2.25}$$

the number of nearest-neighbor spacings in the sequence $x_1, \ldots, x_N$ below the threshold $s$, respectively the number of these spacings between those of the $x_1, \ldots, x_N$ that are located in the interval with diameter $2t$ around $u$. Clearly, $S^{(N)}(s; x) = \lim_{t \to \infty} S_t^{(N)}(s, u; x)$. It is convenient to extend $S^{(N)}(s; \cdot)$ and $S_t^{(N)}(s, u; \cdot)$ to permutation invariant functions on $\mathbb{R}^N$.

**Lemma 2.7.** *For any $N \in \mathbb{N}$ and $t, s > 0$, and $u \in \mathbb{R}$,*

$$\mathbb{E}\big[S_t^{(N)}(s, u; \lambda)\big]$$
$$= \sum_{k=2}^N \frac{(-1)^k}{(k-1)!} \int_{u-t}^{u+t} dr \int_{[0,s]^{k-1}} R_k^{(N)}(r, r+y_2, r+y_3, \ldots, r+y_k) \, dy_2 \cdots dy_k. \tag{2.26}$$

**Sketch of the proof.** We do this only for $t = \infty$. For $k \geq 2$ and $y = (y_1, \ldots, y_k) \in \mathbb{R}^k$, let

$$\chi_{k,s}(y) = \prod_{i,j=1}^k \mathbb{1}\{|y_i - y_j| \leq s\} \quad \text{and}$$

$$Z_{k,s}^{(N)}(y) = \sum_{1 \leq j_1 < \cdots < j_k \leq N} \chi_{k,s}(y_{j_1}, \ldots, y_{j_k}).$$

Elementary combinatorial considerations show that $S^{(N)}(s; x) = \sum_{k=2}^N (-1)^k Z_{k,s}^{(N)}(x)$ for any $x \in W_N$. Multiplying this with the density $P_N$, integrating over $W_N$ and using the permutation symmetry of $\widehat{P}_N = \frac{1}{N!} P_N$ and $Z_{k,s}^{(N)}$ yields

$$\mathbb{E}\big[S^{(N)}(s; \lambda)\big] = \sum_{k=2}^N (-1)^k \int_{W_k} \mathbb{1}\{x_k - x_1 \leq s\} R_k^{(N)}(x) \, d^k x. \tag{2.27}$$

Now an obvious change of variables and the symmetry of $R_k^{(N)}$ yields the assertion for $t = \infty$. □

### 2.7 The orthogonal polynomial method

In this section we briefly describe the most fruitful and most commonly used ansatz for the deeper investigation of the density $P_N$ in (2.3): the *method of orthogonal polynomials*. This technique has been first applied to random matrices by Mehta [Me60] but relies on much older research. For the general theory of orthogonal polynomials see [Sz75]. We follow [De98, Sect. 5] and treat a general orthogonal polynomial ensemble of the form

$$\widehat{P}_N(x) = \frac{1}{N! Z_N} \Delta_N(x)^2 \prod_{i=1}^N e^{-Q(x_i)}, \qquad x = (x_1, \ldots, x_N) \in \mathbb{R}^N, \tag{2.28}$$

where $Q \colon \mathbb{R} \to \mathbb{R}$ is continuous and so large at infinity that all moments of the measure $e^{-Q(x)} \, dx$ are finite. We normalized $\widehat{P}_N$ to a probability density on $\mathbb{R}^N$.

Let $(\pi_j)_{j \in \mathbb{N}_0}$ with $\pi_j(x) = x^j + b_{j-1} x^{j-1} + \cdots + b_1 x + b_0$ be the sequence of orthogonal polynomials for the measure $e^{-Q(x)} \, dx$, i.e.,

$$\int_{\mathbb{R}} \pi_i(x) \pi_j(x) e^{-Q(x)} \, dx = c_i c_j \delta_{ij}, \qquad i, j \in \mathbb{N}_0. \tag{2.29}$$

(In the GUE-case $Q(x) = x^2$, these are the well-known *Hermite polynomials*, scaled such that the leading coefficients are one.) Elementary linear manipulations show that the Vandermonde determinant in (1.3) can be expressed in terms of the same determinant with the monomials $x^j$ replaced by the polynomials $\pi_j(x)$, i.e.,

$$\Delta_N(x) = \det\bigl[(\pi_{j-1}(x_i))_{i,j=1,\ldots,N}\bigr], \qquad x \in \mathbb{R}^N. \tag{2.30}$$

We normalize the $\pi_j$ now to obtain an orthonormal system $(\phi_j)_{j \in \mathbb{N}_0}$ of $L^2(\mathbb{R})$ with respect to the Lebesgue measure: the functions

$$\phi_j(x) = \frac{1}{c_j} e^{-Q(x)/2} \pi_j(x) \tag{2.31}$$

satisfy

$$\int_\mathbb{R} \phi_i(x)\phi_j(x)\, dx = \delta_{ij}, \qquad i,j \in \mathbb{N}_0. \tag{2.32}$$

An important role is played by the kernel $K_N$ defined by

$$K_N(x,y) = \sum_{j=0}^{N-1} \phi_j(x)\phi_j(y), \qquad x,y \in \mathbb{R}. \tag{2.33}$$

The $k$-point correlation function $R_k^{(N)}$ in (2.22) admits the following fundamental determinantal representation.

**Lemma 2.8.** *Fix $N \in \mathbb{N}$ and $x \in \mathbb{R}^N$, then, for any $k \in \{1, \ldots, N\}$,*

$$R_k^{(N)}(x_1, \ldots, x_k) = \det\bigl[(K_N(x_i, x_j))_{i,j=1,\ldots,k}\bigr]. \tag{2.34}$$

*In particular,*

$$R_1^{(N)}(x_1) = K_N(x_1, x_1) \qquad \text{and} \qquad \widehat{P}_N(x) = \frac{1}{N!} \det\bigl[(K_N(x_i, x_j))_{i,j=1,\ldots,N}\bigr]. \tag{2.35}$$

**Sketch of the proof.** Using the determinant multiplication theorem, is easily seen that the density $\widehat{P}_N$ may be written in terms of the functions $\phi_j$ as

$$\begin{aligned}\widehat{P}_N(x) &= \frac{c_0^2 c_1^2 \ldots c_{N-1}^2}{N! Z_N} \det\bigl[(\phi_{j-1}(x_i))_{i,j=1,\ldots,N}\bigr]^2 \\ &= \frac{1}{\widetilde{Z}_N} \det\bigl[(K_N(x_i, x_j))_{i,j=1,\ldots,N}\bigr],\end{aligned} \tag{2.36}$$

where $\widetilde{Z}_N = N! Z_N \prod_{i=0}^{N-1} c_i^{-2}$. Using the special structure of this kernel and some elegant but elementary integration method (see [De98, Lemma 5.27]), one sees that the structure of the density is preserved under successive integration over the coordinates, i.e.,

$$\int_{\mathbb{R}^{N-k}} \det\bigl[(K_N(x_i, x_j))_{i,j=1,\ldots,N}\bigr]\, dx_{k+1} dx_{k+2} \ldots dx_N$$
$$= (N-k)! \det\bigl[(K_N(x_i, x_j))_{i,j=1,\ldots,k}\bigr], \qquad 1 \leq k \leq N. \tag{2.37}$$

In particular, $\widetilde{Z}_N = N!$, and (2.34) holds. □

**Remark 2.9 (Determinantal processes).** Lemma 2.8 offers an important opportunity for far-reaching generalisations. One calls a point process (i.e., a locally finite collection of random points on the real line) a *determinantal process* if its correlation functions are given in the form (2.34), where $K$ is, for some measure $\mu$ on $\mathbb{R}$ having all moments, the kernel of a nonnegative and locally trace class integral operator $L^2(\mathbb{R}, \mu) \to L^2(\mathbb{R}, \mu)$. Because of [De98, Lemma 5.27], correlation functions that are built according to (2.34) form a consistent family of $N$-particle distributions and therefore determine a point process on $\mathbb{R}$. To a certain extent, random matrix calculations only depend on the determinantal structure of the correlation functions are may be used as a starting point for generalisations. ◇

Now let $\lambda = (\lambda_1, \ldots, \lambda_N) \in W_N$ be a random variable with density $P_N = N!\widehat{P}_N$. We now express the probability that a given interval $I$ contains a certain number of $\lambda_i$'s in terms of the operator $\mathcal{K}_N$ on $L^2(I)$ with kernel $K_N(x,y)$.

**Lemma 2.10.** *For any $N \in \mathbb{N}$, any $m \in \{0, \ldots, N\}$, and any interval $I \subset \mathbb{R}$,*

$$\mathbb{P}\big(\#\{i \leq N \colon \lambda_i \in I\} = m\big) = \frac{(-1)^m}{m!}\Big(\frac{\mathrm{d}}{\mathrm{d}\gamma}\Big)^m \det\big[(\mathrm{Id} - \gamma\mathcal{K}_N)|_{L^2(I)}\big]\Big|_{\gamma=1}, \quad (2.38)$$

*where* Id *denotes the identical operator in $L^2(I)$.*

**Sketch of the proof.** From Lemma 2.6 and (2.34) we have

$$\begin{aligned}&\mathbb{P}\big(\#\{i \leq N \colon \lambda_i \in I\} = m\big) \\ &= \frac{(-1)^m}{m!} \sum_{k=m}^{N} \frac{(-1)^k}{(k-m)!} \int_{I^k} \det\big[(K_N(x_i, x_j))_{i,j=1,\ldots,k}\big]\, \mathrm{d}^k x.\end{aligned} \quad (2.39)$$

On the other hand, for any $\gamma \in \mathbb{R}$, by a classical formula for trace class operators, see [RS7580, Vol. IV, Sect. 17]

$$\det\big[(\mathrm{Id} - \gamma\mathcal{K}_N)|_{L^2(I)}\big] = \sum_{k=0}^{N} \frac{(-\gamma)^k}{k!} \int_{I^k} \det\big[(K_N(x_i, x_j))_{i,j=1,\ldots,k}\big]\, \mathrm{d}^k x. \quad (2.40)$$

Now differentiate $m$ times with respect to $\gamma$ and put $\gamma = 1$ to arrive at (2.38). □

### 2.8 Spacings in the bulk of the spectrum, and the sine kernel

In this section, we explain the limiting spacing statistics in the bulk of the spectrum of a random unitary-invariant $(N \times N)$ matrix in the limit $N \to \infty$. We specialize to the matrix distribution in (2.1) with $F$ as in (2.2) and $Q(x) = x^{2j}$ for some $j \in \mathbb{N}$. This has the technical advantage of a perfect-scaling property of the eigenvalues: as was pointed out in Remark 2.5(ii), the correct rescaling of the eigenvalues is $\widetilde{\lambda}_i^{(N)} = N^{-\frac{1}{2j}} \lambda_i^{(N)}$. In order to ease the notation, we shall consider $\widetilde{\lambda}^{(N)}$ instead of $\lambda^{(N)}$. Note that the distribution of $\widetilde{\lambda}^{(N)}$ is the orthogonal polynomial ensemble in (2.28) with $Q(x) = Nx^{2j}$, and we shall stick to that choice of $Q$ from now.

Let $\psi \colon \mathbb{R} \to [0, \infty)$ denote the density of the equilibrium measure (i.e., the unique minimizer) for the functional in (2.19) with $\gamma_{2j} = 1$. According to the semicircle law, the rescaled eigenvalues $\widetilde{\lambda}_i^{(N)}$ lie asymptotically in the support of $\psi$, which is the rescaled bulk of the spectrum. In particular, the spacings between neighboring rescaled eigenvalues should be of order $\frac{1}{N}$, and hence the spacings of the unscaled eigenvalues are of order $N^{\frac{1}{2j}-1}$.

We fix a threshold $s > 0$ and a point $u \in \mathrm{supp}(\psi)^\circ$ in the bulk of the rescaled spectrum and want to describe the number of spacings $\leq \frac{s}{N}$ of the rescaled eigenvalues in a vicinity of $u$. Let $(t_N)_{N \in \mathbb{N}}$ be a sequence in $(0, \infty)$ with $t_N \to 0$ as $N \to \infty$. The main object of our interest is the expected value of $S_{t_N}^{(N)}(\frac{s}{N}, u; \widetilde{\lambda}^{(N)})$, the number of spacings $\leq \frac{s}{N}$ in the sequence $\widetilde{\lambda}^{(N)}$ in a $t_N$-interval around $u$; see (2.25). We expect that this number is comparable to $t_N N$, and we want to find the asymptotic dependence on $s$ and $u$.

We continue to follow [De98, Sect. 5] and stay in the framework of Section 2.7, keeping all assumptions and all notation, and specializing to $Q(x) = Nx^{2j}$. We indicate the $N$-dependence of the weight function $Q(x) = Nx^{2j}$ by writing $K_N^{(N)}$ for the kernel $K_N$ defined in (2.33) and (2.31). Abbreviate

$$\kappa_N(u) = K_N^{(N)}(u, u). \tag{2.41}$$

We write $\widetilde{R}_1^{(N)}$ for the 1-point correlation function with respect to the ensemble in (2.28) with $Q(x) = Nx^{2j}$; hence $\widetilde{R}_1^{(N)}(u)\,\mathrm{d}u$ is the density of $\frac{1}{N}$ times the number of rescaled eigenvalues in $\mathrm{d}u$ (see below (2.22)). From (2.35) we have $\kappa_N(u) = \widetilde{R}_1^{(N)}(u)$. Hence, the asymptotics of $\kappa_N(u)$ can be guessed from the semi-circle law: we should have $\kappa_N(u) = \widetilde{R}_1^{(N)}(u) \approx N\psi(u)$. (In the GUE-case $j = 1$, we have $|u| < \sqrt{2}$ and $\kappa_N(u) \approx N\frac{1}{\pi}\sqrt{2 - u^2}$.) We shall adapt the scaling of the expected number of spacings to the spot $u$ where they are registered by using the scaling $\frac{1}{\kappa_N(u)}$ instead of $\frac{1}{N}$. This will turn out to make the value of the scaling limit independent of $u$.

We use now Lemmas 2.7 and 2.8 and an elementary change of the integration variables to find the expectation of the number of rescaled eigenvalue spacings as follows.

$$\frac{1}{\kappa_N(u)} \frac{1}{2t_N} \mathbb{E}\Big[S_{t_N}^{(N)}\Big(\frac{s}{\kappa_N(u)}, u; \widetilde{\lambda}^{(N)}\Big)\Big] = \sum_{k=2}^{N} \frac{(-1)^k}{(k-1)!} \frac{1}{2t_N} \int_{u-t_N}^{u+t_N} \mathrm{d}r$$
$$\int_{[0,s]^{k-1}} \det\Big[\Big(\frac{1}{\kappa_N(u)} K_N^{(N)}\Big(r + \frac{y_i}{\kappa_N(u)}, r + \frac{y_j}{\kappa_N(u)}\Big)\Big)_{i,j=1,\dots,k}\Big]\Big|_{y_1=0} \mathrm{d}y_2 \dots \mathrm{d}y_k. \tag{2.42}$$

Hence, we need the convergence of the rescaled kernel in the determinant on the right hand side. This is provided in the following theorem. The well-known *sine*

*kernel* is defined by

$$\mathbb{S}(x - y) = \frac{\sin(\pi(x - y))}{\pi(x - y)} = \frac{\sin(\pi x)\sin'(\pi y) - \sin'(\pi x)\sin(\pi y)}{\pi(x - y)}, \qquad x, y \in \mathbb{R}. \tag{2.43}$$

**Proposition 2.11 (Bulk asymptotics for $K_N^{(N)}$).** *Fix $j \in \mathbb{N}$ and $Q(x) = Nx^{2j}$. Let $K_N^{(N)}$ be as in (2.33) with the functions $\phi_j$ defined in (2.31) such that (2.32) holds. Denote by $\psi \colon \mathbb{R} \to [0, \infty)$ the equilibrium measure of the functional in (2.19) with $\gamma_{2j} = 1$. For $u \in \mathrm{supp}(\psi)^\circ$, abbreviate $\kappa_N(u) = K_N^{(N)}(u, u)$. Then, uniformly on compact subsets in $u \in \mathrm{supp}(\psi)^\circ$ and $x, y \in \mathbb{R}$,*

$$\lim_{N \to \infty} \frac{1}{\kappa_N(u)} K_N^{(N)}\!\left(u + \frac{x}{\kappa_N(u)}, u + \frac{y}{\kappa_N(u)}\right) = \mathbb{S}(x - y). \tag{2.44}$$

For a rough outline of the proof using Riemann-Hilbert theory, see Section 2.10 below.

**Remark 2.12.** (i) The asymptotics in (2.44) in the GUE case $j = 1$, where the orthogonal polynomials are the Hermite polynomials, are called the *Plancherel-Rotach asymptotics*.
 (ii) Note that the limit in (2.44) is independent of $u$, as long as $u$ is in the interior of the support of the equilibrium measure, i.e., as long as we consider the bulk of the spectrum. See Proposition 2.15 for the edge asymptotics.
 (iii) The asymptotics in Proposition 2.11 are universal in the sense that they do not depend on the weight function $Q(x)$, at least within the class $Q(x) = x^{2j}$ we consider here (after rescaling). The case of a polynomial $Q(x) = x^{2j} + \mathcal{O}(x^{2j-1})$ is asymptotically the same, but the proof is technically more involved.
 (iv) The proof of Proposition 2.11 is in [De98, Ch. 8], based on [KS99]. The first proof, even for more general functions $Q$, is in [PS97]. See also [D99] and [BI99] for related results. The main tool for deriving (2.44) (and many asymptotic assertions about orthogonal polynomials) are the *Riemann-Hilbert theory* and the *Deift-Zhou steepest decent method*.
 (v) Analogous results for weight functions of Laguerre type (recall (2.6)) for $\beta = 2$ have been derived using adaptations of the methods mentioned in (iv). The best available result seems to be in [Va05], where weight functions of the form $\mu(\mathrm{d}x_i) = x_i^\alpha e^{-Q(x_i)}\,\mathrm{d}x_i$ are considered with $\alpha > -1$, and $Q$ is an even polynomial with positive leading coefficient. The cases $\beta = 1$ and $\beta = 4$ are considered in [DGKV05].
 (vi) The orthogonal and symplectic cases (i.e., $\beta = 1$ and $\beta = 4$) for Hermite-type weight functions $\mu(\mathrm{d}x_i) = e^{-Q(x_i)}\,\mathrm{d}x_i$ with $Q$ a polynomial have also been carried out recently [DG05a].
 (vii) Using the well-known *Christoffel-Darboux formula*

$$\sum_{j=0}^{N-1} q_j(x) q_j(y) = \frac{c_N}{c_{N-1}} \frac{q_N(x) q_{N-1}(y) - q_N(y) q_{N-1}(x)}{x - y}, \qquad x, y \in \mathbb{R}, \tag{2.45}$$

(where $q_j = \pi_j/c_j$; see (2.31)), one can express the kernel $K_N$ defined in (2.33) in terms of just two of the orthogonal polynomials. Note the formal analogy between the right hand sides of (2.45) and (2.43).    $\diamond$

Now we formulate the main assertion about the limiting eigenvalue spacing for random unitary-invariant matrices. Denote by $\mathcal{K}_{\sin}$ the integral operator whose kernel is the sine kernel in (2.43).

**Theorem 2.13 (Limiting eigenvalue spacing).** *Let $M_N$ be a random unitary-invariant matrix of the form (2.1) with $F$ as in (2.2) and $Q(x) = x^{2j}$ for some $j \in \mathbb{N}$. Let $\widetilde{\lambda}^{(N)} = N^{-\frac{1}{2j}}(\lambda_1^{(N)}, \ldots, \lambda_N^{(N)})$ be the vector of scaled eigenvalues of $M_N$. Denote by $\psi \colon \mathbb{R} \to [0, \infty)$ the equilibrium measure of the functional in (2.19) with $\gamma_{2j} = 1$. Fix $u \in \mathrm{supp}(\psi)^\circ$ and $s > 0$ and let $(t_N)_{N \in \mathbb{N}}$ be a sequence in $(0, \infty)$ with $t_N \to 0$. Recall the definition (2.25) of the spacing number. Then*

$$\lim_{N \to \infty} \frac{1}{N\psi(u)} \frac{1}{2t_N} \mathbb{E}\left[ S_{t_N}^{(N)}\left(\frac{s}{N\psi(u)}, u; \widetilde{\lambda}^{(N)}\right) \right] = \int_0^s p(v)\,dv, \qquad (2.46)$$

*where*

$$p(v) = \frac{d^2}{dv^2} \det\left[ (\mathrm{Id} - \mathcal{K}_{\sin})|_{L^2([v, \infty))} \right], \qquad v \geq 0, \qquad (2.47)$$

*is the density of the* Gaudin *distribution.*

**Sketch of the proof.** In (2.42), replace the normalized $r$-integral by the integral over the delta-measure on $u$ and use Proposition 2.11 to obtain

$$\text{left hand side of (2.46)}$$
$$= \sum_{k=2}^{\infty} \frac{(-1)^k}{(k-1)!} \int_{[0,s]^{k-1}} \det\left[ (\mathbb{S}(y_i - y_j))_{i,j=1,\ldots,k} \right]\Big|_{y_1=0} dy_2 \cdots dy_k. \qquad (2.48)$$

On the other hand, note that

$$\int_0^s p(v)\,dv = 1 + \frac{d}{ds} \det\left[ (\mathrm{Id} - \mathcal{K}_{\sin})|_{L^2([s, \infty))} \right]$$

$$= 1 - \frac{d}{d\varepsilon}\Big|_{\varepsilon=0} \det\left[ (\mathrm{Id} - \mathcal{K}_{\sin})|_{L^2([\varepsilon, s])} \right]$$

$$= 1 - \frac{d}{d\varepsilon}\Big|_{\varepsilon=0} \left[ 1 + \sum_{k=1}^{\infty} \frac{(-1)^k}{k!} \int_{[\varepsilon, s]^k} \det\left[ (\mathbb{S}(y_i - y_j))_{i,j=1,\ldots,k} \right] d^k y \right]$$

$$= \text{right hand side of (2.48)},$$
$$(2.49)$$

as an application of the product differentiation rule shows.    $\square$

**Remark 2.14.** (i) It is instructive to compare the asymptotic spacing distribution of the rescaled eigenvalues of a large random matrix (which have a mutual repellence) to the one of $N$ independently on the interval $[0, 1]$ randomly and uniformly distributed points (where no interaction appears).

The latter can be realized as a standard conditional Poisson process, given that there are precisely $N$ Poisson points in $[0,1]$. The asymptotic spacing density for the latter is just $v \mapsto e^{-v}$, and the former is $v \mapsto p(v)$ as in Theorem 2.13. Note that the asymptotics of $p(v)$ for $v \downarrow 0$ and the one for $v \to \infty$ are both smaller than the one of $e^{-v}$. Indeed, it is known that $p(v) \approx v^\beta$ as $v \downarrow 0$ if in (2.3) the term $\Delta_N(x)^2$ is replaced by $\Delta_N(x)^\beta$ and, furthermore, $p(v) \approx e^{-v^2}$ as $v \to \infty$; see [DIZ97] and [De98, Sect. 8.2].

(ii) Another variant of the assertion in (2.46) is about the number of pairs of rescaled, not necessarily neighboring, eigenvalues whose difference is in a fixed interval $(a,b)$:

$$\lim_{N \to \infty} \frac{1}{N} \mathbb{E}\Big[\#\big\{(i,j) \in \{1,\ldots,N\}^2 \colon a < \widetilde{\lambda}_i^{(N)} - \widetilde{\lambda}_j^{(N)} < b\big\}\Big]$$
$$= \int_a^b \Big[1 - \Big(\frac{\sin(\pi u)}{\pi u}\Big)^2\Big] \,\mathrm{d}u + \begin{cases} 1 & \text{if } 0 \in (a,b), \\ 0 & \text{otherwise.} \end{cases} \tag{2.50}$$

The last term accounts for the pairs $i = j$.

(iii) Proposition 2.11 and Theorem 2.13 are extended to a large class of Wigner matrices in [Jo01a], more precisely to the class of random Hermitian matrices of the form $W + aV$, where $W$ is a Wigner matrix as in Section 2.2, $a > 0$ and $V$ is a standard GUE-matrix. The entries of $W$ are not assumed to have a symmetric distribution, but the expected value is supposed to be zero, the variance is fixed, and the $(6+\varepsilon)$-th moments for any $\varepsilon > 0$ are supposed to be uniformly bounded. This result shows universality of the limiting spacing distribution in a large class of Wigner matrices. The identification of the distribution of the eigenvalues of $W + aV$ uses the interpretation of the eigenvalue process of $(W + aV)_{a \geq 0}$ as a process of non-colliding Brownian motions as in [Dy62b], see Section 4.1 below.

(iv) After appropriate asymptotic centering and normalization, the distribution of the individual eigenvalues for GUE in the bulk of the spectrum is asymptotically Gaussian. Indeed, for $i_N = (a+o(1))N$ with $a \in (-\sqrt{2}, \sqrt{2})$ (i.e., $a$ is in the interior of the support of the semicircle law $\mu_*$ in (2.10)), the correct scaling of the $i_N$-th eigenvalue is

$$X_{i_N}^{(N)} = \frac{\lambda_{i_N}^{(N)} - t\sqrt{2N}}{\big(\frac{\log N}{(1-2t^2)N}\big)^{1/2}},$$

where $t$ is determined by $\mu_*((-\infty, t]) = a$. One main result of [Gus04] is that $X_{i_N}^{(N)}$ is asymptotically standard normal as $N \to \infty$. Also joint distributions of several bulk eigenvalues in this scaling are considered in [Gus04]. In particular, it turns out that $\lambda_i^{(N)}$ and $\lambda_j^{(N)}$ are asymptotically independent if $|i-j|$ is of the order $N$, but not if $|i-j| = o(N)$. ◇

## 2.9 The edge of the spectrum, and the Tracy-Widom distribution

In this section we explain the limiting scaled distribution of the largest eigenvalue, $\lambda_N^{(N)}$, of an $(N \times N)$ GUE-matrix, i.e., we specialize to $j = 1$. Let $\lambda^{(N)} = (\lambda_1^{(N)}, \ldots, \lambda_N^{(N)}) \in W_N$ be the vector of the eigenvalues. According to Lemma 2.1, its distribution is the orthogonal polynomial ensemble in (2.28) with $Q(x) = x^2$. Hence, the distribution of the vector of rescaled eigenvalues, $N^{-1/2}\lambda^{(N)}$, is that ensemble with $Q(x) = Nx^2$. The event $\{\lambda_N^{(N)} \leq t\}$ is, for any $t \in \mathbb{R}$, identical to the event that no eigenvalue falls into the interval $(t, \infty)$. Hence we may apply Lemma 2.10 for $I = (t, \infty)$ and $m = 0$. In order to obtain an interesting limit as $N \to \infty$, we already know from the semicircle law that $t$ should be chosen as $t = \sqrt{2N} + \mathcal{O}(N^\alpha)$ for some $\alpha < \frac{1}{2}$. It will turn out that $\alpha = -\frac{1}{6}$ is the correct choice.

As in the preceding section, we denote by $K_N^{(N)}$ the kernel $K_N$ defined in (2.33) for the choice $Q(x) = Nx^2$, with the functions $\phi_j$ defined in (2.31) such that (2.32) holds. Using Lemma 2.6 for $m = 0$ and (2.34), we see, after an elementary change of measure, that

$$\mathbb{P}\bigl(\lambda_N^{(N)} \leq \sqrt{2N} + s(\sqrt{2}N^{\frac{1}{6}})^{-1}\bigr)$$
$$= \mathbb{P}\Bigl(\lambda_1^{(N)}, \ldots, \lambda_N^{(N)} \notin \bigl(\sqrt{2N} + s(\sqrt{2}N^{\frac{1}{6}})^{-1}, \infty\bigr)\Bigr)$$
$$= \sum_{k=0}^N \frac{(-1)^k}{k!} \int_{[s,\infty)^k} \det\Bigl[\Bigl(\frac{1}{\sqrt{2}N^{\frac{2}{3}}} K_N^{(N)}\Bigl(\sqrt{2} + \frac{x_i}{\sqrt{2}N^{\frac{2}{3}}}, \sqrt{2} + \frac{x_j}{\sqrt{2}N^{\frac{2}{3}}}\Bigr)\Bigr)_{i,j=1,\ldots,k}\Bigr] \mathrm{d}^k x. \quad (2.51)$$

In order to obtain an interesting limit, one needs to show that the integrand on the right hand side of (2.51) converges. This is provided in the following theorem. By Ai: $\mathbb{R} \to \mathbb{R}$ we denote the *Airy function*, the unique solution to the differential equation $f''(x) = xf(x)$ on $\mathbb{R}$ with asymptotics $f(x) \sim (4\pi\sqrt{x})^{1/2} e^{-\frac{2}{3}x^{3/2}}$ as $x \to \infty$. The corresponding kernel, the *Airy kernel*, is given by

$$K_{\mathrm{Ai}}(x,y) = \frac{\mathrm{Ai}(x)\mathrm{Ai}'(y) - \mathrm{Ai}'(x)\mathrm{Ai}(y)}{x - y} = \int_0^\infty \mathrm{Ai}(x+u)\mathrm{Ai}(y+u)\,\mathrm{d}u, \qquad x, y \in \mathbb{R}. \quad (2.52)$$

Note the formal analogy to (2.43) and (2.45).

**Proposition 2.15 (Edge asymptotics for $K_N$).** *Uniformly in $x, y \in \mathbb{R}$ on compacts,*

$$\lim_{N\to\infty} \frac{1}{\sqrt{2}N^{\frac{2}{3}}} K_N^{(N)}\Bigl(\sqrt{2} + \frac{x}{\sqrt{2}N^{\frac{2}{3}}}, \sqrt{2} + \frac{y}{\sqrt{2}N^{\frac{2}{3}}}\Bigr) = K_{\mathrm{Ai}}(x,y). \quad (2.53)$$

**Remark 2.16.** (i) Note that the kernel $K_N^{(N)}$ scales with $N^{-\frac{2}{3}}$ at the edge of the spectrum, i.e., in $\pm\sqrt{2}$, while it scales with $\frac{1}{N}$ in the interior of the support of the equilibrium measure, $(-\sqrt{2}, \sqrt{2})$ (see Proposition 2.11).
 (ii) The Airy kernel already appeared in [BB91] in a related connection. Proofs of Proposition 2.15 were found independently by Tracy and Widom [TW93a] and Forrester [Fo93].
 (iii) For an extension of Proposition 2.15 to the weight function $Q(x) = x^{2j}$ for some $j \in \mathbb{N}$, see [De98, Sec. 7.6], e.g. The real and symplectic cases (i.e., $\beta = 1$ and $\beta = 4$) have also been recently carried out [DG05b].

(iv) Analogous results for weight functions of Laguerre type (recall (2.6) and Remark 2.12(v)) for $\beta = 1$ and $\beta = 4$ are derived in [DGKV05]. Both boundaries, the 'hard' edge at zero and the 'soft ' one at the other end, are considered. ◇

Next, we formulate the asymptotics for the edge of the spectrum, i.e., the largest (resp. smallest) eigenvalues. Let $q\colon \mathbb{R} \to \mathbb{R}$ be the solution[5] [HML80] of the Painlevé II differential equation

$$q''(x) = xq(x) + 2q(x)^3 \qquad (2.54)$$

with asymptotics $q(x) \sim \mathrm{Ai}(x)$ as $x \to \infty$. It is uniquely determined by requiring that $q(x) > 0$ for any $x < 0$, and it has asymptotics $q(x) \sim \sqrt{|x|/2}$ as $x \to -\infty$. Furthermore, $q'(x) < 0$ for any $x \in \mathbb{R}$.

Define a distribution function $F_2\colon \mathbb{R} \to [0,1]$ by

$$F_2(s) = \exp\Big\{-\int_s^\infty (x-s)q^2(x)\,\mathrm{d}x\Big\}, \qquad s \in \mathbb{R}. \qquad (2.55)$$

This is the distribution of the by now famous *GUE Tracy-Widom distribution*; its importance is clear from the following.

**Theorem 2.17 (Limiting distribution of the largest eigenvalue, [TW94a]).** *Let $M_N$ be a random Hermitian matrix of size $N$ from GUE, and let $\lambda_N^{(N)}$ be the largest eigenvalue of $M_N$. Then, for any $s \in \mathbb{R}$,*

$$\lim_{N\to\infty} \mathbb{P}\Big(\big(\lambda_N^{(N)} - \sqrt{2N}\big)\sqrt{2}N^{1/6} \leq s\Big) = F_2(s). \qquad (2.56)$$

**Proof.** Using (2.51) and Proposition 2.15, we obtain

$$\lim_{N\to\infty} \mathbb{P}\Big(\big(\lambda_N^{(N)} - \sqrt{2N}\big)\sqrt{2}N^{1/6} \leq s\Big)$$
$$= \sum_{k=0}^\infty \frac{(-1)^k}{k!} \int_{[s,\infty)^k} \det\Big[(K_{\mathrm{Ai}}(x_i, x_j))_{i,j=1,\ldots,k}\Big]\,\mathrm{d}^k x \qquad (2.57)$$
$$= \det\Big[(\mathrm{Id} - \mathcal{K}_{\mathrm{Ai}})\big|_{L^2([s,\infty))}\Big],$$

where $\mathcal{K}_{\mathrm{Ai}}$ is the operator on $L^2([s,\infty))$ with kernel $K_{\mathrm{Ai}}$. The relation to the Painlevé equation is derived in [TW94a] using a combination of techniques from operator theory and ordinary differential equations. □

**Remark 2.18.** (i) The great value of Theorem 2.17 is the characterization of the limit on the left hand side of (2.57) in terms of some basic ordinary differential equation, in this case the Painlevé II equation. Analogous relations between the Gaudin distribution $p$ in (2.47) and the Painlevé V equation were derived in [JMMS80].

---
[5]The function $u \equiv -q$ is also a solution of (2.54), which is sometimes called the *Hastings-Mac Leod solution*.

(ii) There are analogous statements for GOE and GSE [TW96]. The limiting distributions are modifications of the GUE Tracy-Widom distribution. Indeed, for $\beta = 1$ and $\beta = 4$, respectively (recall (2.5)), $F_2$ is replaced by

$$F_1(s) = \exp\Big\{-\frac{1}{2}\int_s^\infty \big[q(x) + (x-s)q^2(x)\big]\,\mathrm{d}x\Big\} = \sqrt{F_2(s)}\,e^{-\frac{1}{2}\int_s^\infty q(x)\,\mathrm{d}x},$$

$$F_4(s) = \sqrt{F_2(s)}\,\frac{1}{2}\big[e^{\frac{1}{2}\int_s^\infty q(x)\,\mathrm{d}x} + e^{-\frac{1}{2}\int_s^\infty q(x)\,\mathrm{d}x}\big].$$

(2.58)

(iii) The expectation of a random variable with distribution function $F_2$ is negative and has the value of approximately $-1.7711$.

(iv) In [TW94a], also the joint distribution of the first $m$ top eigenvalues is treated; they admit an analogous limit theorem. The starting point for the proof is Lemma 2.6 and (2.34).

(v) Theorem 2.17 is generalized to a large class of Wigner matrices in [So99]. It is assumed there that the entries have a symmetric distribution with all moments finite such that the asymptotics for high moments are bounded by those of the Gaussian. The proof is a variant of the method of moments (see the first proof of Theorem 2.3). The main point is that the expected trace of high powers (appropriately coupled with the matrix size) of the random matrix is bounded by a certain asymptotics, which is essentially the same as for GUE. Since the expected trace of high moments depends on the matrix entries only via the moments, which are the same within the class considered, the result then follows from a comparison to the known asymptotics for GUE.

(vi) If the index $i_N$ is a bit away from the edge $N$, then the $i_N$-th largest eigenvalue scales to some Gaussian law. Indeed, if $i_N = N - k_N$ with $k_N \to \infty$, but $k_N/N \to 0$, then the correct scaling is

$$X_{i_N}^{(N)} = \frac{\lambda_{i_N}^{(N)} - \sqrt{2N}\big(1 - \big(\frac{3\pi k_N}{4\sqrt{2}N}\big)^{2/3}\big)}{\big((12\pi)^{-2/3}\frac{\log k_N}{N^{1/3}\,k_N^{2/3}}\big)^{1/2}},$$

and one main result of [Gus04] is that $X_{i_N}^{(N)}$ is asymptotically standard normal. Also joint distributions of several eigenvalues in this scaling are considered in [Gus04]. In particular, it turns out that $\lambda_{i_N}^{(N)}$ and $\lambda_{j_N}^{(N)}$ (provided that $N - i_N$ and $N - j_N$ are $o(N)$) are asymptotically correlated if $|i_N - j_N| = o(N - i_N)$. ◇

## 2.10 Some elements of Riemann-Hilbert theory

Apparently, the most powerful technical tool for deriving limiting assertions about orthogonal polynomial ensembles is the *Riemann-Hilbert (RH) theory*. This theory dates back to the 19th century and was originally introduced for the study of monodromy questions in ordinary differential equations, and has been applied to a host of models in analysis. Applications to orthogonal polynomials were first developed in [FIK90], and this method was first combined with a steepest-decent method in [DZ93]. Since then, a lot of deep results on random matrix theory and related models could be established using a combination of the two methods. The reformulation in terms of RH theory leaves the intuition of orthogonal polynomial ensembles behind, but creates a new framework, in which a new intuition arises and new technical tools become applicable which are suitable to deal with the difficulties stemming from the great number of zeros of the polynomials. For a recent general survey on Riemann-Hilbert theory, see [It03]; for a thorough exposition of the adaptation and application of this theory to orthogonal polynomials, see the lectures [De98], and [Ku03], [D01] and [BI99].

In this section, we give a rough indication of how to use Riemann-Hilbert theory for scaling limits for orthogonal polynomials, in particular we outline some elements of the proof of Proposition 2.11. We follow [De98]. Let us start with the definition of a Riemann-Hilbert problem in a situation specialized to our purposes, omitting all technical issues.

Let $\Sigma$ be a finite union of the images of smooth, oriented curves in $\mathbb{C}$, and suppose there is a smooth function $v$ (called the *jump matrix*) on $\Sigma$ with values in the set of complex regular $(2 \times 2)$-matrices. We say a matrix-valued function $Y$ on $\mathbb{C} \setminus \Sigma$ *solves the Riemann-Hilbert (RH) problem* $(\Sigma, v)$ if

$$
\begin{array}{lrll}
(i) & Y & \text{is analytic in } \mathbb{C} \setminus \Sigma, & \\
(ii) & Y_+(x) = Y_-(x)v(x), & x \in \Sigma, & \quad (2.59) \\
(iii) & Y(z) = I + \mathcal{O}(\tfrac{1}{z}) & \text{as } z \to \infty,
\end{array}
$$

where $I$ is the $(2 \times 2)$-identity matrix, and $Y_+(x)$ and $Y_-(x)$ are the limiting boundary values of $Y$ in $x \in \Sigma$ coming from the positive and negative side of $\Sigma$, respectively.[6]

The main connection with orthogonal polynomials is in Proposition 2.19 below. Assume that $\mu(\mathrm{d}x) = w(x)\,\mathrm{d}x$ is a positive measure on $\mathbb{R}$ having all moments and a sufficiently regular density $w$, and let $(\pi_n)_{n \in \mathbb{N}_0}$ be the sequence of orthogonal polynomials for the $L^2$-inner product with weight $w$, such that the degree of $\pi_n$ is $n$ and the highest coefficient one. Hence, for some $k_n > 0$,

$$\int_\mathbb{R} \pi_n(x)\pi_m(x)\,\mu(\mathrm{d}x) = \frac{1}{k_n^2}\delta_{n,m}, \qquad n,m \in \mathbb{N}_0. \quad (2.60)$$

Recall the *Cauchy transform*,

$$\mathcal{C}f(z) = \int_\mathbb{R} \frac{f(x)}{x-z}\frac{\mathrm{d}x}{2\pi\mathrm{i}}, \qquad z \in \mathbb{C} \setminus \mathbb{R}, f \in H^1(\mathbb{R}). \quad (2.61)$$

Here we think of $\mathbb{R}$ as of an oriented curve from $-\infty$ to $\infty$, parametrized by the identity map. Note that $\mathcal{C}(f)_+ - \mathcal{C}(f)_- = f$ on $\mathbb{R}$.

---

[6]The definition of $Y_+(x)$ and $Y_-(x)$ and the sense in which (ii) is to be understood have to be explained rigorously, and (ii) is required outside the intersections of the curves only, but we neglect these issues here. The general notion involves $(k \times k)$-matrices for some $k \in \mathbb{N}$ instead of $(2 \times 2)$-matrices.

**Proposition 2.19 (RH problem for orthogonal polynomials, [FIK90], [FIK91]).** *Fix $n \in \mathbb{N}$ and consider the jump matrix $v(x) = \begin{pmatrix} 1 & w(x) \\ 0 & 1 \end{pmatrix}$ for $x \in \mathbb{R}$. Then*

$$Y^{(n)}(z) = \begin{pmatrix} \pi_n(z) & \mathcal{C}(\pi_n w)(z) \\ -2\pi\mathrm{i}\, k_{n-1}^2 \pi_{n-1}(z) & -2\pi\mathrm{i}\, k_{n-1}^2 \pi_{n-1}(z)\mathcal{C}(\pi_{n-1} w)(z) \end{pmatrix}, \qquad z \in \mathbb{C}\setminus\mathbb{R}. \tag{2.62}$$

*is the unique solution of the RH problem[7]*

$$\begin{array}{rll}
(i) & Y^{(n)} \text{ is analytic in } \mathbb{C} \setminus \mathbb{R}, & \\
(ii) & Y_+^{(n)}(x) = Y_-^{(n)}(x) v(x), & x \in \mathbb{R}, \\
(iii) & Y^{(n)}(z) = (I + \mathcal{O}(\tfrac{1}{z})) \begin{pmatrix} z^n & 0 \\ 0 & z^{-n} \end{pmatrix} & \text{as } z \to \infty.
\end{array} \tag{2.63}$$

The main advantage of the characterisation of the orthogonal polynomials in terms of a RH problem is that it provides a technical frame in which the difficulties stemming from the oscillations of the polynomials close to their zeros can be resolved.

Now we specialize to $w(x) = e^{-NQ(x)}$ with $Q(x) = x^{2j}$ for some $j \in \mathbb{N}$, recall Remark 2.5(ii) and Section 2.8. We now write $\pi_n^{(N)}$ instead of $\pi_n$ for the orthogonal polynomials. We shall (extremely briefly) indicate how the asymptotics of the $N$-th orthogonal polynomial $\pi_N^{(N)}$ can be deduced from RH theory, building on Proposition 2.19.

The first main step is a transformation of (2.63) which absorbs the exponential term of the jump matrix into an inverse exponential term in the solution of the new RH problem.[8] For doing this, we need to use some information about the variational formula in (2.19) with $\gamma_{2j} = 1$. Recall the Euler-Lagrange equations in (2.20) for the equilibrium measure $\mu_*(\mathrm{d}x) = \psi(x)\,\mathrm{d}x$, and put

$$g(z) = \int_\mathbb{R} \log(z-x)\psi(x)\,\mathrm{d}x, \qquad z \in \mathbb{C}\setminus\mathbb{R}. \tag{2.64}$$

The intuitive idea behind the choice of $g$ is the fact that, if $x_1^*, \ldots, x_N^* \in \mathbb{R}$ denote the zeros of $\pi_N^{(N)}$ and $\mu_N$ their empirical measure, then we can write

$$\pi_N^{(N)}(x) = \prod_{i=1}^N (x - x_i^*) = e^{N\int_\mathbb{R} \log(x-y)\,\mu_N(\mathrm{d}y)}; \tag{2.65}$$

compare also to (1.4). Since the asymptotic statistics of the zeros and of the ensemble particles are very close to each other, we should have $\pi_N^{(N)} \approx e^{Ng}$, and $e^{Ng}$ will indeed turn out to be the main term in the expansion.

Consider the transformed jump matrix

$$v^{(1)}(x) = \begin{pmatrix} e^{N[g_-(x)-g_+(x)]} & e^{N[g_-(x)+g_+(x)-Q(x)+l]} \\ 0 & e^{N[g_+(x)-g_-(x)]} \end{pmatrix}, \qquad x \in \mathbb{R}. \tag{2.66}$$

---

[7]Note that (2.63) is not a standard RH problem, compare (iii).
[8]This step is analogous to the exponential change of measure in large deviation theory, which absorbs the main exponential factor in the probability with respect to a transformed measure under which the deviation behavior becomes typical. However, because of the great amount of zeros of $\pi_N^{(N)}$, the exponential term is rather subtle here.

Then the unique solution, $m^{(1)}$, of the RH problem $(\mathbb{R}, v^{(1)})$ can easily be calculated from $Y^{(n)}$ in Proposition 2.19; its $(1,1)$-entry is $\pi_N^{(N)} e^{-Ng}$. This means that the leading (exponential) term has been isolated in the transformed RH problem $(\mathbb{R}, v^{(1)})$. It turns out that, outside the support of the equilibrium measure, $v^{(1)}(x)$ is exponentially close to the identity matrix, and inside this support we have

$$v^{(1)}(x) = \begin{pmatrix} e^{-2\pi i\, N\mu_*([x,\infty))} & 1 \\ 0 & e^{2\pi i\, N\mu_*([x,\infty))} \end{pmatrix}$$
$$= \begin{pmatrix} 1 & 0 \\ e^{2N\varphi_-(x)} & 1 \end{pmatrix} \begin{pmatrix} 0 & 1 \\ -1 & 0 \end{pmatrix} \begin{pmatrix} 1 & 0 \\ e^{2N\varphi_+(x)} & 1 \end{pmatrix}, \tag{2.67}$$

where $\varphi(z) = 2\pi i \int_a^z \psi(t)\,dt$ for $z \in \mathbb{C} \setminus [-a, a]$, where we chose $a > 0$ such that $\mathrm{supp}(\mu_*) = [-a, a]$. We know from Remark 2.5(ii) that $t \mapsto \psi(t)(a^2 - t^2)^{-\frac{1}{2}}$ is analytic in $\mathbb{C}$, and hence $\int_a^z \psi(t)\,dt$ depends on the integration curve from $a$ to $z$: any two curves lead to a difference by an integer multiple of $2\pi i$. Hence, $z \mapsto e^{n\varphi(z)}$ is well-defined and analytic in $\mathbb{C} \setminus [-a, a]$ and therefore this is true for its boundary functions on $(-a, a)$, $\varphi_+$ and $\varphi_-$.

The next main step is a deformation of $(\mathbb{R}, v^{(1)})$, which isolates the second term in the expansion of $\pi_N^{(N)}$, which is of fourth-root order and hence much more subtle. Indeed, the decomposition in the second line of (2.67) gives rise to a deformation into a RH problem $(\Sigma, v^{(2)})$, where $\Sigma$ is the union of the real line and two curves from $-a$ to $a$ in the upper and lower half plane, respectively, and $v^{(2)}$ is some suitable jump matrix on $\Sigma$. It is relatively easy to prove that, in $L^2$-sense, as $N \to \infty$, we have $v^{(2)} \to v_\infty$ on $\Sigma$ with $v_\infty = \begin{pmatrix} 0 & 1 \\ -1 & 0 \end{pmatrix}$ on $[-a, a]$. Hence, the unique solution, $m^{(2)}$, of the problem $(\Sigma, v^{(2)})$ should converge towards the unique solution, $m_\infty$, of the RH problem $([-a, a], v_\infty)$. This is true, but relatively hard to prove, in particular on $\mathrm{supp}(\mu_*)$ and here in particular close to the boundaries $\pm a$. It is easy to compute that

$$m_\infty = \frac{1}{2i}\begin{pmatrix} i(\beta + \beta^{-1}) & \beta - \beta^{-1} \\ \beta^{-1} - \beta & i(\beta + \beta^{-1}) \end{pmatrix}, \qquad \text{where } \beta(z) = \left(\frac{z-a}{z+a}\right)^{\frac{1}{4}}. \tag{2.68}$$

Computing $m^{(2)}$, re-substituting $m^{(1)}$ and $m_\infty$, and considering the $(1,1)$-entry, we obtain therefore the asymptotics of $\pi_N^{(N)}$ outside the critical points $\pm a$:

$$\pi_N^{(N)}(z) = \begin{cases} \frac{1}{2}\left[\left(\frac{z-a}{z+a}\right)^{\frac{1}{4}} + \left(\frac{z+a}{z-a}\right)^{\frac{1}{4}} + o(1)\right] e^{Ng(z)} \\ \qquad \text{if } z \in \mathbb{C} \setminus \mathrm{supp}(\mu_*), \\ \frac{1}{2}\left[\left|\frac{z-a}{z+a}\right|^{\frac{1}{4}} \cos\left(N\pi\mu_*([z,a]) + \frac{\pi}{4}\right) \\ \quad + \left|\frac{z+a}{z-a}\right|^{\frac{1}{4}} \cos\left(N\pi\mu_*([z,a]) - \frac{\pi}{4}\right) + o(1)\right] e^{N\int_{-a}^a \log|z-x|\,\mu_*(dx)} \\ \qquad \text{if } z \in \mathrm{supp}(\mu_*)^\circ. \end{cases} \tag{2.69}$$

This explains how to derive the Plancherel-Rotach asymptotics for the orthogonal polynomials for the weight function $w(x) = e^{-Nx^{2j}}$. Note that the error terms in (2.69) are locally uniform outside neighborhoods of $\pm a$. Exploiting the Christoffel-Darboux formula in (2.45), one obtains the statement of Proposition 2.11.

In order to obtain the asymptotics of Proposition 2.15, i.e., the asymptotics of $\pi_N^{(N)}(z)$ for $z$ close to $\pm a$, one uses an appropriate deformation into a suitable RH problems involving the Airy function, see [De98, Sect. 7.6], e.g.

**2.11 Random matrices and the Riemann zeta function**

Excitingly, it turned out in the early seventies that the spacings of the zeros of the Riemann zeta function show a close relation to those of the eigenvalues of certain random matrices. The famous *Riemann zeta function* is defined on $\{\Re(s) > 1\}$ as

$$\zeta(s) = \sum_{n=1}^{\infty} n^{-s} = \prod_{p \, \text{prime}} (1 - p^{-s})^{-1}. \tag{2.70}$$

Riemann showed in 1859 that $\zeta$ can be meromorphically continued to the whole complex plane, and that the functional equation $\Gamma(s/2)\zeta(s)\sqrt{\pi} = \pi^s \Gamma(\frac{1}{2}(1-s))\zeta(1-s)$ holds. This continuation has simple zeros at the negative even integers and a simple pole at 1, which is the only singularity. Furthermore, there are infinitely many zeros in the so-called critical strip $\{0 < \Re(s) < 1\}$, and none of them is real. These zeros are called the non-trivial zeros; they are located symmetrically around the real axis and around the line $\{\Re(s) = \frac{1}{2}\}$, the critical line. Denote them by $\rho_n = \beta_n + \mathrm{i}\,\gamma_n$ with $\gamma_{-1} < 0 < \gamma_1 \leq \gamma_2 \leq \dots$. The famous *Riemann Hypothesis* conjectures that $\beta_n = \frac{1}{2}$ for every $n$, i.e., every non-trivial zero lies on the critical line $\{\Re(s) = \frac{1}{2}\}$. This is one of the most famous open problems in mathematics and has far reaching connections to other branches of mathematics. Many rigorous results in analytic number theory are conditional on the truth of the Riemann Hypothesis. There is extensive evidence for it being true, as many partial rigorous results and computer simulations have shown. See [Ed74] and [Ti86] for much more on the Riemann zeta function.

It is known that the number $\pi(x)$ of prime numbers $\leq x$ behaves asymptotically as $\pi(x) = \mathrm{Li}(x) + \mathcal{O}(x^\Theta \log x)$ as $x \to \infty$, where $\mathrm{Li}(x)$ is the principal value of $\int_0^x \frac{1}{\log u}\,\mathrm{d}u$, which is asymptotic to $\frac{x}{\log x}$, and $\Theta = \sup_{n \in \mathbb{N}} \beta_n$. Hence, the Riemann Hypothesis is equivalent to a precise asymptotic statement about the prime number distribution.

More interestingly from the viewpoint of orthogonal polynomial ensembles, the Riemann Hypothesis has also much to do with the vertical distribution of the Riemann zeros. Let $N(T)$ be the number of zeros in the critical strip up to height $T$, counted according to multiplicity. It is known that $N(T) = \frac{T}{2\pi} \log \frac{T}{2\pi e} + \mathcal{O}(\log T)$ as $T \to \infty$. In the pioneering work [Mo73], vertical spacings of the Riemann zeros are considered. Denote by

$$R_T(a,b) = \frac{1}{N(T)} \#\Big\{(n,m) \in \mathbb{N}^2 \colon \gamma_n, \gamma_m \leq T, a \leq \frac{\gamma_n - \gamma_m}{2\pi} \log \frac{T}{2\pi} \leq b\Big\}, \qquad a < b, \tag{2.71}$$

the number of pairs of rescaled critical Riemann zeros whose difference lies between $a$ and $b$. Then it was proved in [Mo73], assuming the Riemann Hypothesis,

that

$$\lim_{T\to\infty} R_T(a,b) = \int_a^b \left[1 - \left(\frac{\sin(\pi u)}{\pi u}\right)^2\right] du + \begin{cases} 1 & \text{if } 0 \in (a,b), \\ 0 & \text{otherwise.} \end{cases} \quad (2.72)$$

The last term accounts for the pairs $m = n$. Note the close similarity to (2.50). Calculating millions of zeros, [Od87] confirms this asymptotics with an extraordinary accuracy.

The *Lindelöf Hypothesis* is the conjecture that $\zeta(\frac{1}{2}+\mathrm{i}\,t) = \mathcal{O}(t^\varepsilon)$ as $t \to \infty$ for any $\varepsilon > 0$. The $(2k)$-th moment of the modulus of the Riemann zeta function,

$$I_k(T) = \frac{1}{T} \int_0^T \left|\zeta(\tfrac{1}{2} + \mathrm{i}\,t)\right|^{2k} dt, \quad (2.73)$$

was originally studied in an attempt to prove the Lindelöf Hypothesis, which is equivalent to $I_k(T) = \mathcal{O}(t^\varepsilon)$ as $T \to \infty$ for any $\varepsilon > 0$ and any $k$. The latter statement has been proved for $k = 1$ and $k = 2$. Based on random matrix calculations, [KS00] conjectured that

$$I_k(T) \sim \frac{G^2(k+1)}{G(2k+1)}\, a(k)(\log T)^{k^2}, \qquad k \in \{\Re(s) > -\tfrac{1}{2}\}, \quad (2.74)$$

where $G$ is the *Barnes G-function*, and

$$a(k) = \prod_{p\,\text{prime}} \left(1 - \frac{1}{p}\right)^{k^2} \sum_{m \in \mathbb{N}_0} \left(\frac{\Gamma(m+k)}{m!\,\Gamma(k)}\right)^2 p^{-m}. \quad (2.75)$$

This so-called *Keating-Snaith Conjecture* was derived by an asymptotic calculation for the Fourier transform of the logarithm of the characteristic polynomial of a random matrix from the Circular Unitary Ensemble introduced in Section 2.4. This conjecture is one of the rare (non-rigorous, however) progresses in recent decades in the understanding of the Riemann zeros.

## 3. Random growth processes

In this section we consider certain classes of random growth processes which turned out in the late 1990es to be closely connected to certain orthogonal polynomial ensembles. There is a number of physically motivated random growth processes which model growing surfaces under influences of randomly occurring events (like nucleation events) that locally increase a substrate, but have far-reaching correlations on a long-time run. In one space dimension, for these kinds of growth processes, limiting phenomena are conjectured that have morally some features of random matrices in common, like the fluctuation behavior of power-order $1/3$ (instead of the order $1/2$ in the central limit theorem and related phenomena) and the universality of certain rescaled quantities. Recently some of these models could be analysed rigorously, after exciting discoveries of surprising relations to orthogonal polynomial ensembles had been made.

Random growth models may be defined in any dimension, and two and three dimensional models are of high interest. However, the high-dimensional cases seem mathematically intractable yet, such that we restrict to one-dimensional[9] models in this text. General physics references on growing surfaces are the monographs [BS95] and [Me98]; see also [KS92]. Much background is also provided in [P03] and [Fe04b]. Recent surveys on some growth models that have been solved in recent years by methods analogous to those used in random matrix theory are [Jo01c] and [Ba03].

After a short description of one basic model that cannot be handled rigorously yet in Section 3.1, we shall treat basically only two models: the corner-growth model introduced in Section 3.2 and the PNG model introduced in Section 3.6. The main results on these two models are presented in Sections 3.3 and Section 3.4, respectively in Sections 3.6 and 3.7. The famous and much-studied problem of the longest increasing subsequence of a random permutation is touched in Section 3.5, since it is instrumental for the PNG model (and also important on its own). Furthermore, in Section 3.8, we mention the Plancherel measure as an technically important toy model that links combinatorics and orthogonal polynomials.

## 3.1 The Eden-Richardson model

A fundamental model for random growth is the so-called *Eden-Richardson model*, which is defined as follows. The model is a random process $(A(t))_{t\geq 0}$ of subsets of $\mathbb{Z}^2$ such that $A(t) \subset A(s)$ for any $t < s$. At time $t = 0$, the set $A_0$ is equal to $\{0\}$, the origin in $\mathbb{Z}^2$. We call a site $(i,j) \in \mathbb{Z}^2 \setminus A(t)$ *active* at time $t$ if some neighbor of $(i,j)$ belongs to $A(t)$. As soon as $(i,j)$ is active, a random waiting time $w(i,j)$ starts running, and after this time has elapsed, $(i,j)$ is added to the set process as well. The waiting times $w(i,j)$, $(i,j) \in \mathbb{Z}^2$, are assumed to be independent and identically distributed $(0,\infty)$-valued random variables. They can be discrete or continuous. In the case of $\mathbb{N}$-valued waiting times, we consider the discrete-time process $(A(t))_{t\in\mathbb{N}_0}$ instead of $(A(t))_{t\geq 0}$. If and only if the distribution of the waiting times is exponential, respectively geometric, the process $(A(t))_{t\geq 0}$, respectively $(A(t))_{t\in\mathbb{N}_0}$, enjoys the Markov property: in the discrete-time case, at each time unit any active site chooses independently with a fixed probability if it immediately belongs to the set process or not. In this special case, the model is called the Eden-Richardson model. The Markov property is not present for any other distribution.

Actually, the Eden-Richardson model is equivalent to what probabilists call *last-passage percolation*, which we will explain more closely in Remark 3.1 below.

The natural question is about the asymptotic behavior of the set $A(t)$ for large $t$. It is not so difficult to conjecture that there should be a law of large numbers be valid, i.e., there should be a deterministic set $A \subset \mathbb{R}^2$ such that $\frac{1}{t}A(t) \to A$ as $t \to \infty$. A proof of this fact can be derived using the subadditive ergodic theorem [Ke86], which considers the Markovian case. However, an identification of the limiting set $A$ and closer descriptions of $A$ for general waiting time distributions seem out of reach. In physics literature, it is conjectured that the fluctuations be of order $t^{1/3}$. It is rather hard to analyze Eden's model mathematically rigorously. Reasons for that are that $A(t)$ may and does have holes and that the growth proceeds in any direction. No technique has yet been found to attack the asymptotics of the fluctuations rigorously. This is why we do not spend time on the Eden model, but immediately turn to some simpler variant which has been successfully treated.

---
[9]Taking into acount the time-evolution, they are sometimes also called (1+1)-dimensional.

### 3.2 The corner-growth model

An important simpler variant of Eden's model is known as the *corner growth model*. This is a growth model on $\mathbb{N}_0^2$ instead of $\mathbb{Z}^2$, and growth is possible only in corners. At time zero, $A(0)$ is the union of the $x$-axis $\mathbb{N}_0 \times \{0\}$ and the $y$-axis $\{0\} \times \mathbb{N}_0$. Points in $\mathbb{N}^2 \setminus A(t)$ are called *active* at time $t$ if their left and their lower neighbors both belong to $A(t)$. As soon as a point $(i,j)$ is active, its individual waiting time $w(i,j)$ starts running, and after it elapses $(i,j)$ is added to the set. This defines a random process $(A(t))_{t \geq 0}$ of growing subsets of $\mathbb{N}_0^2$. Again, if the waiting times are $\mathbb{N}$-valued, we consider $(A(t))_{t \in \mathbb{N}_0}$, and the Markov property is present only for the two above mentioned waiting time distributions: the exponential, respectively the geometric, distributions.

It is convenient to identify every point $(i,j)$ with the square $[i-\frac{1}{2}, i+\frac{1}{2}) \times [j-\frac{1}{2}, j+\frac{1}{2})$ and to regard $A(t)$ as a subset of $[\frac{1}{2}, \infty)^2$. The process $(A(t))_{t \geq 0}$ consists of an infinite number of growing columns, of which almost all are of zero height and which are ordered non-increasingly in height. One can view these columns as a vector of runners who proceed like independent random walkers, making a unit step after an individual independent waiting time, subject to the rule that the $(i+1)$-st runner is stopped by the $i$-th runner as long as they are on the same level. Note that this is a *suppression* mechanism, not a *conditioning* mechanism. A realization of $A(t)$ is as follows (the active sites are marked by '×').

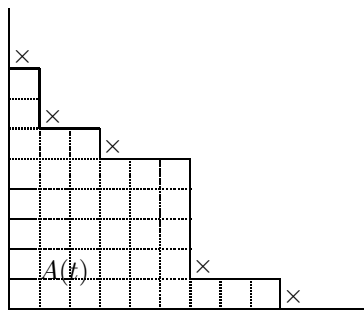

Denote by $G(M, N)$ the first time $t$ at which the point $(M, N)$ belongs to $A(t)$. Obviously, $G(M, N)$ depends on the variables $w(i,j)$ with $i \leq M$ and $j \leq N$ only, and the recurrence relation

$$G(M, N) = w(M, N) + \max\{G(M-1, N), G(M, N-1)\}, \qquad M, N \in \mathbb{N}, \quad (3.1)$$

is satisfied. From this relation, it is straightforward to derive the formula

$$G(M, N) = \max_{\pi \in \Pi(M,N)} \sum_{(i,j) \in \pi} w(i,j), \qquad M, N \in \mathbb{N}, \tag{3.2}$$

where $\Pi(M, N)$ is the set of nearest-neigbor paths with $M+N-2$ steps, starting from $(1, 1)$, ending at $(M, N)$ and having only steps $(0, 1)$ or $(1, 0)$. An example of an array of waiting times with $M = 4$ and $N = 6$ is

| 2 | 0 | 2 | 5 | 2 | 5 |
|---|---|---|---|---|---|
| 2 | 6 | 3 | 4 | 1 | 6 |
| 3 | 0 | 0 | 0 | 3 | 1 |
| 2 | 7 | 3 | 4 | 1 | 0 |

We ordered the rows from the bottom to the top. Four maximal paths from the left lower corner to the right upper corner are depicted.

Much of the interest in the corner-growth model stems from the fact that it has a couple of connections to other well-known models and admits several alternate descriptions:

**Remark 3.1 (Directed last-passage percolation).** Switch to the dual lattice of $\mathbb{Z}^2$ (where the role of sites and bonds are interchanged) and interpret $w(i, j)$ as the travel time along the bond $(i, j)$. Then $\sum_{(i,j) \in \pi} w(i, j)$ is the travel time along the path $\pi \in \Pi(M, N)$ from $(0, 0)$ to $(M, N)$. Hence, $G(M, N)$ is the maximal travel time from $(0, 0)$ to $(M, N)$ along right/up paths, and this model is called directed last-passage percolation. The term 'directed' refers to the use of right/up paths. If *all* nearest-neighbor paths are admissable, then the model is called last-passage percolation; this is equivalent to the Eden-Richardson model of Section 3.1.

Switching the signs of $w(i, j)$ and ignoring that $-w(i, j)$ is negative, we see that $-G(M, N)$ is the *minimal* travel time (now with passage 'times' $-w(i, j)$) from $(0, 0)$ to $(M, N)$, which is the well-known model of *first-passage percolation*. An interpretation is as follows. If at the origin there is the source of a fluid, whose floating time along the bond $(i, j)$ is $-w(i, j)$, then the set $\widetilde{A}(t) = \{(M, N) \colon -G(M, N) \leq t\}$ is the set of bonds that are wet by time $t$. ◇

**Remark 3.2 (Totally asymmetric exclusion process).** The boundary of the set $A(t) \subset [\frac{1}{2}, \infty)^2$ is a curve that begins with infinitely many vertical line segments of unit length, proceeds with finitely many horizontal and vertical line segments of unit length, and ends with infinitely many horizontal line segments of unit length. If a square is added to $A(t)$, then a vertical/horizontal pair of lines is changed into a horizontal/vertical pair. If we replace vertical lines by a '1' and horizontal lines by a '0' and determine the index that refers to the main diagonal of $\mathbb{R}^2$ as 0, then we can think of the corner growth model as of a particle process $(x_k(t))_{k \in \mathbb{Z}} \in \{0, 1\}^{\mathbb{Z}}$ where $x_k(t) = 1$ means that one particle is present at site $k$ at time $t$. In the case of geometric waiting time distribution, the dynamics of this process is as follows. At each discrete time unit, every particle independently moves to the right neighboring site with a fixed probability, provided this site is vacant. Otherwise, it does not move. These are the dynamics of the so-called totally asymmetric exclusion process in discrete time. The event $\{G(M, N) = t\}$ is the event that the particle that was initially at site $1 - N$ has moved $M$ steps by time $t$. There is an analogous representation in continuous time for the exponential waiting time distribution. ◇

**Remark 3.3 (Directed polymers in random environment).** Let $(S_n)_{n \in \mathbb{N}_0}$ be a simple random walk on $\mathbb{Z}$, then the process $(n, S_n)_{n \in \mathbb{N}_0}$ is interpreted as a directed polymer in $\mathbb{Z}^2$. Let $(v(i, j))_{i \in \mathbb{N}_0, j \in \mathbb{Z}}$ be an i.i.d. field of real random variables. Every monomer $(n, S_n)$ receives the weight $\beta v(n, S_n)$, where $\beta > 0$ is interpreted as the inverse of the temperature. This induces a probability measure on $N$-step paths given by

$$Q_{N,\beta}(S_0, \ldots, S_N) = \frac{1}{Z_{N,\beta}} \exp\Big\{-\beta \sum_{n=0}^{N} v(n, S_n)\Big\}.$$

In the zero-temperature limit $\beta \to \infty$, the measure $Q_{N,\beta}$ is concentrated on those paths $(S_0, \ldots, S_N)$ which minimize $\sum_{n=0}^{N} v(n, S_n)$. This is the analog of the corner-growth model with switched signs of the random variables; compare to (3.2). It is believed that the directed polymer at positive, sufficiently small, temperature essentially exhibits the same large-$N$ behavior as the zero-temperature limit, but this conjecture is largely unproven. An account on the recent research on directed polymers in random environment is in [CSY04]. ◇

**Remark 3.4 (Tandem queues).** At time zero, there is an infinite number of customers in the first queue, and there is an infinite number of other queues, which are initially empty and have to be passed by every customer one after another. The first customer in any queue (if present) is served after a random waiting time, which has the distribution of the waiting times in the corner growth model, and he or she proceeds to the end of the next queue. Then, at every time $t$, the height of the $i$-th column of the set $A(t)$ is equal to the number of customers which have passed or reached the $i$-th queue. A general and systematic discussion of the relation between tandem queues and orthogonal polynomial ensembles appears in [OC03]. ◇

A systematic study of the random variable on the right side of (3.2) and its asymptotics towards Brownian analogs is in [Ba01]; see also [GTW01], [OY02], [BJ02] and [Do03]. In fact, for $N$ fixed and under appropriate moment conditions, in the limit $M \to \infty$, this random variable (after proper centering and rescaling) converges in distribution towards

$$L(N) = \max_{1 \geq t_1 \geq \cdots \geq t_{N-1} \geq 0} \Big[W_1(1) - W_1(t_1) + W_2(t_1) - W_2(t_2) \pm \ldots \\ \pm W_N(t_{N-1}) - W_N(0)\Big], \quad (3.3)$$

where $W_1, \ldots, W_N$ are $N$ independent standard Brownian motions on $\mathbb{R}$ starting at the origin. Using Donsker's invariance principle, this may be explained as follows. Assume that $\mathbb{E}[w(1,1)] = 0$ and $\mathbb{E}[w(1,1)^2] = 1$. The first upstep of a path $\pi$ in (3.2) may be expected in the $(t_{N-1}M)$-th step, the second in the $(t_{N-2}M)$-th step and so on, where we later optimize on $1 \geq t_1 \geq \cdots \geq t_{N-1} \geq 0$. The partial sums of $w(i, t_{i-1}M), \ldots, w(i, t_iM - 1)$ approach the distribution of $\sqrt{M}[W_i(t_{i-1}) - W_i(t_i)]$ for $i = 1, \ldots, N$. Hence, we have that $M^{-1/2}G(M, N) \Longrightarrow L(N)$ as $M \to \infty$, in the case of centered and normalized random variables $w(i, j)$.

A rather beautiful fact [Ba01], [GTW01] is that $L(N)$ is in distribution equal to the largest eigenvalue of a GUE matrix, $\lambda_N^{(N)}$. (For generalisations of this fact to Brownian motion in the fundamental chamber associated with a finite Coxeter group see [BBO05].) Recall from Theorem 2.17 that we may approximate $\lambda_N^{(N)} \approx \sqrt{2N} + (\sqrt{2}N^{\frac{1}{6}})^{-1}F_2$ for $N$ large. Combining the limits for $M \to \infty$ and $N \to \infty$, one is lead to the appealing conjecture (still assuming that $\mathbb{E}[w(1,1)] = 0$ and $\mathbb{E}[w(1,1)^2] = 1$)

$$\sqrt{2}N^{\frac{1}{6}}\Big(\frac{G(M,N)}{\sqrt{M}} - \sqrt{2N}\Big) \Longrightarrow F_2, \qquad M, N \to \infty. \quad (3.4)$$

This assertion has indeed been proven independently in [BM05] and [BS05], under the additional assumption that $M = o(N^a)$ for $a < \frac{3}{14}$. The main tool is a classical strong approximation of random walks by Brownian motion, which works so well that $M$ can diverge together with $N$ at some speed. However, the most interesting case is where $M$ and $N$ are of the same order, and this case is open yet in general. For the two special cases of the geometric and the exponential distribution, (3.4) has been proven for $M \approx \text{const.} \times N$. Our next two sections are devoted to a description of this result.

### 3.3 Johansson's identification of the distribution

In his beautiful work [Jo00a], Kurt Johansson deeply investigated the corner-growth model for both particular waiting-time distributions, the geometric and the exponential distribution. He identified the distribution of $G(M, N)$ in terms of the distribution of the largest particle of the *Laguerre ensemble* (see (2.6)) in the exponential case, and of the *Meixner ensemble* (its discrete analog) in the geometric case.

**Proposition 3.5 (Distribution of $G(M, N)$, [Jo00a]).** *Let $G(M, N)$ be defined as in (3.2), and let the $w(i, j)$ be i.i.d. geometrically distributed with parameter $q \in (0, 1)$, i.e., $w(i, j) = k \in \mathbb{N}$ with probability $(1 - q)q^k$. Then, for any $M, N \in \mathbb{N}$ with $M \geq N$, and for any $t \in \mathbb{N}$,*

$$\mathbb{P}(G(M,N) \leq t) = \frac{1}{Z_{M,N}} \sum_{x_1,\ldots,x_N=1}^{t+N-1} \Delta_N(x)^2 \prod_{i=1}^N \Big[\binom{x_i + M - N}{x_i} q^{x_i}\Big]. \quad (3.5)$$

**Remark 3.6.**   (i) The right hand side of (3.5) is the probability that the largest particle in the *Meixner ensemble* on $\mathbb{N}^N$ with parameters $q$ and $M - N$ is smaller than $t + N$.

(ii) There is an extension of Proposition 3.5 to the case where the parameter of the geometric distribution of $w(i,j)$ is of the form $a_i b_j$ for certain numbers $a_i, b_j \in (0,1)$, see [Jo01c, Sect. 2].
(iii) An analogous formula holds for the case of exponentially distributed waiting times, and the corresponding ensemble is the Laguerre ensemble (Gamma-distribution in place of the negative Binomial distribution), see (2.6). This formula is derived in [Jo00a] using an elementary limiting procedure which produces the exponential distribution from the geometric one. It is remarkable that no direct proof is known yet. Distributions other than the exponential or geometric one have not yet been successfully treated. ◇

**Sketch of the proof of Proposition 3.5.** The proof in [Jo00a] relies on various combinatoric tools, which have been useful in various parts of mathematics for decades. A general reference is [Sa91].

A *generalized permutation* is an array of two rows with integer entries such that the row of the pairs is non-decreasingly ordered in lexicographical sense. An example is

$$\sigma = \begin{pmatrix} 1 & 1 & 1 & 1 & 1 & 2 & 2 & 2 & 2 & 3 & 4 & 4 \\ 1 & 1 & 3 & 3 & 3 & 1 & 1 & 1 & 3 & 3 & 2 & 3 \end{pmatrix}; \quad (3.6)$$

the entries of the first and second line are taken from $\{1,2,3,4\}$ and $\{1,2,3\}$, respectively. A longest increasing subsequence of the second row has the length 8; it consists of all the '1' and the last three '3'. Also the first two ones and all the threes form a longest increasing subsequence.

**Lemma 3.7 (Matrices and generalized permutations).** *For any $M, N, k \in \mathbb{N}$, the following procedure defines a one-to-one map between the set of $(M \times N)$-matrices $(W(i,j))_{i \leq M, j \leq N}$ with positive integer entries and total sum $\sum_{i \leq M, j \leq N} W(i,j)$ equal to $k$, and the set of generalized permutations of length $k$ whose first row has entries in $\{1, \ldots, M\}$ and whose second row has entries in $\{1, \ldots, N\}$: Repeat every pair $(i,j) \in \{1, \ldots, M\} \times \{1, \ldots, N\}$ precisely $W(i,j)$ times, and list all pairs in lexicographical order. By this procedure, the quantity $\max_{\pi \in \Pi(M,N)} \sum_{(i,j) \in \pi} W(i,j)$ is mapped onto the length of the longest non-decreasing subsequence of the second row.*

As an example for $M = 4$, $N = 3$, the matrix

$$W = \begin{pmatrix} 0 & 1 & 1 \\ 0 & 0 & 1 \\ 3 & 0 & 1 \\ 2 & 0 & 3 \end{pmatrix} \quad (3.7)$$

is mapped onto the generalized permutation $\sigma$ in (3.6). (In order to appeal to the orientation of the corner growth model, we ordered the rows of $W$ from the bottom to the top, contrary to the order one is used to from linear algebra.) The two paths linking the coordinates $(1,1), (2,1), (2,3), (4,3)$ and $(1,1), (1,3), (4,3)$, respectively, are maximal paths in (3.2); they correspond to the longest increasing subsequences mentioned below (3.6).

**Remark 3.8.** (i) For the application of Lemma 3.7 for $W(i,j) = w(i,j)$ geometrically distributed random variables, it is of crucial importance that this distribution induces a uniform distribution on the set of $(M \times N)$-matrices with fixed sum of the entries.

(ii) Obviously, Lemma 3.7 works *a priori* only for integer-valued matrices. ◇

The next step is a famous bijection between generalized permutations and Young tableaux. A *semi-standard Young tableau*[10] is a finite array of rows, nonincreasing in lengths, having integer entries which are nondecreasing along the rows and strictly increasing along the columns. The *shape* of the tableau, $\lambda = (\lambda_i)_i$, is the vector of lengths of the rows. In particular, $\lambda_1$ is the length of the longest row of the tableau, and $\sum_i \lambda_i$ is the total number of entries. An example of a semi-standard Young tableau with shape $\lambda = (10, 8, 8, 3, 1)$ and entries in $\{1, \ldots, 6\}$ is as follows.

| 1 | 1 | 2 | 2 | 3 | 3 | 3 | 4 | 4 | 6 |
|---|---|---|---|---|---|---|---|---|---|
| 2 | 2 | 3 | 4 | 4 | 4 | 5 | 5 |   |   |
| 3 | 3 | 5 | 5 | 5 | 5 | 6 | 6 |   |   |
| 4 | 5 | 6 |   |   |   |   |   |   |   |
| 6 |   |   |   |   |   |   |   |   |   |

**Lemma 3.9 (Robinson-Schensted-Knuth (RSK) correspondence, [K70]).** *For any $M, N, k \in \mathbb{N}$, there is a bijection between the set of generalized permutations of length $k$ whose first row has entries in $\{1, \ldots, M\}$ and whose second row has entries in $\{1, \ldots, N\}$, and the set of pairs of semi-standard Young tableaux of the same shape with total number of entries equal to $k$, such that the entries of the first Young tableau are taken from $\{1, \ldots, M\}$ and the ones of the second from $\{1, \ldots, N\}$. This bijection maps the length of the longest non-decreasing subsequence of the second row of the permutation onto the length of the first row of the tableau, $\lambda_1$.*

The algorithm was introduced in [Sc61] for permutations (it is a variant of the well-known *patience sorting algorithm*) and was extended to generalized permutations in [K70].

Sofar, the distribution of $G(M, N)$ has been reformulated in terms of the length of the first row of pairs of semi-standard Young tableaux. The next and final tool is a combinatorial formula for the number of Young tableaux.

**Lemma 3.10 (Number of semi-standard Young tableaux).** *The number of semi-standard Young tableaux of shape $\lambda$ and elements in $\{1, \ldots, N\}$ is equal to*
$$\prod_{1 \leq i < j \leq N} \frac{\lambda_i - \lambda_j + j - i}{j - i}.$$

---

[10]For the notions of (standard) Young tableaux and Young diagrams, see Section 3.8 below.

The reader easily recognizes that the combinatorial formula in Lemma 3.10 is the kernel of the formula in (3.5). Putting together the tools listed sofar, one easily arrives at (3.5). □

**Remark 3.11.** An alternate characterization and derivation of the distribution of $G(M, N)$ is given in [Jo02a, Sect. 2.4] in terms of the *Krawtchouk ensemble*,

$$\mathrm{Kr}_{M,n,q}(x) = \frac{1}{Z_{M,n,q}} \Delta_M(x)^2 \prod_{i=1}^{M} \left[ \binom{n}{x_i} q^{x_i}(1-q)^{n-x_i} \right], \qquad x \in \{0, \ldots, n\}^M \cap W_M. \tag{3.8}$$

There a family of random non-colliding one-dimensional nearest-neighbor processes is analyzed, which is a discrete analog of the multilayer PNG-droplet model in Section 3.7 below. The joint distribution of this cascade of processes is identified in terms of the the Krawtchouk ensemble, and the marginal distribution of the rightmost process is identified in terms of $G(M, N)$. This implies that

$$\mathbb{P}(G(M, N) \leq t) = \sum_{x \in \{0, \ldots, t+M-1\}^M} \mathrm{Kr}_{M, t+N+M-1, q}(x), \tag{3.9}$$

i.e., $G(M, N)$ is characterized in terms of the largest particle of the Krawtchouk ensemble. ◇

### 3.4 Asymptotics for the Markovian corner-growth model

Having arrived at the description in (3.5), the machinery of statistical mechanics and orthogonal polynomials can be applied. The outcome is the following.

**Theorem 3.12 (Asymptotics for the corner-growth model, [Jo00a]).** *Consider the model of Proposition 3.5. Then, for any $\gamma \geq 1$,*

$$(i) \qquad \lim_{N \to \infty} \frac{1}{N} \mathbb{E}\big[G(\lfloor \gamma N \rfloor, N)\big] = \frac{(1 + \sqrt{q\gamma})^2}{1 - q} - 1 \equiv f(\gamma, q), \tag{3.10}$$

$$(ii) \qquad \lim_{N \to \infty} \mathbb{P}\Big( \frac{G(\lfloor \gamma N \rfloor, N) - N f(\gamma, q)}{\sigma(\gamma, q) N^{1/3}} \leq s \Big) = F_2(s), \qquad s \in \mathbb{R}, \tag{3.11}$$

*where $F_2$ is the distribution function of the GUE Tracy-Widom distribution introduced in* (2.55), *and $\sigma(\gamma, q)$ is some explicit function.*

**Remark 3.13.**  (i) In Theorem 3.12 the weak law of large numbers $\lim_{t \to \infty} \frac{1}{t} A(t) = A$ is contained with

$$A = \{(x, y) \in [0, \infty)^2 \colon y + 2\sqrt{qxy} + x \leq 1 - q\}.$$

A qualitative picture of $A$ is as follows.

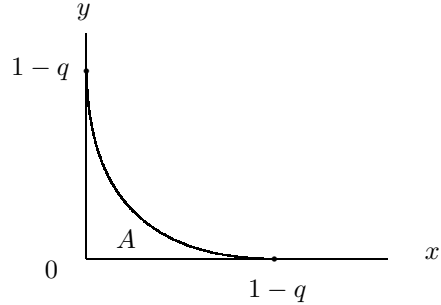

(ii) The analogous result for exponentially distributed waiting times is also contained in [Jo00a].
(iii) Theorem 3.12(ii) is equivalent to (3.4) with $M = \lfloor \gamma N \rfloor$ (recall that the $w(i,j)$ are centered and normalized there, but not in Theorem 3.12).
(iv) Theorem 3.12 is the first and yet only rigorous result of behavior of the type 'fluctuations $\approx$ mean$^{\frac{1}{3}}$' for random growth models of the Eden type.
(v) In [Jo00a] also some large-deviation estimates are proven, i.e., the probabilities of the events $\{G(\lfloor \gamma N \rfloor, N) \leq N(f(\gamma, q) - \varepsilon)\}$ and of $\{G(\lfloor \gamma N \rfloor, N) \geq N(f(\gamma, q) + \varepsilon)\}$ are asymptotically estimated. The former decays on the scale $e^{\mathcal{O}(N^2)}$, while the latter decays on the scale $e^{\mathcal{O}(N)}$.
(vi) The proof of Theorem 3.12 can also be built on the Krawtchouk representation in (3.9) by using asymptotics for the Krawtchouk kernel; see [Jo02a]. ◇

**Sketch of the proof of Theorem 3.12.** The structure of this proof is analogous to the proof of Theorem 2.17. The right hand side of (3.5) may be written in terms of the Meixner kernel

$$K_{\text{Me}}^{(N)}(x,y) = \sum_{j=0}^{N-1} m_j(x) m_j(y) \left[ w_L^{(q)}(x) w_L^{(q)}(y) \right]^{\frac{1}{2}}, \qquad L = M - N + 1, \quad (3.12)$$

where $m_j(x) = \kappa_j x^j + \mathcal{O}(x^{j-1})$ are the orthonormal polynomials with respect to the discrete weight $w_L^{(q)}(x) = \binom{x+L-1}{x} q^x$, $x \in \mathbb{N}$. (Both the polynomials $m_j$ and the kernel $K_{\text{Me}}^{(N)}$ also depend on $L$ and $q$.) Indeed, computations similar to those of Section 2.7 imply that

right hand side of (3.5)

$$= \sum_{k=0}^{N} \frac{(-1)^k}{k!} \sum_{h \in \{t+N, t+N+1, \ldots\}^k} \det \left[ \left( K_{\text{Me}}^{(N)}(h_i, h_j) \right)_{i,j=1,\ldots,k} \right]. \quad (3.13)$$

The Meixner kernel satisfies the scaling limit

$$\lim_{N \to \infty} \sigma N^{\frac{1}{3}} K_{\text{Me}}^{(N)} \left( (f+1)N + \sigma N^{\frac{1}{3}} x, (f+1)N + \sigma N^{\frac{1}{3}} y \right) = K_{\text{Ai}}(x,y), \quad (3.14)$$

where $K_{\text{Ai}}$ is the Airy kernel in (2.53), and $f = f(\gamma, q)$ and $\sigma = \sigma(\gamma, q)$ are as in the theorem. Now the remainder of the proof is analogous to the proof of Theorem 2.17. □

### 3.5 Longest increasing subsequences of random permutations

Another problem that has been recognized to be strongly related to random growth processes is the problem of the length of the longest increasing subsequence of a random permutation. Let $\mathfrak{S}_N$ denote the set of permutations of $1, \ldots, N$, and let $\sigma$ be a random variable that is uniformly distributed on $\mathfrak{S}_N$, i.e., a random permutation. The length of the longest increasing subsequence of $\sigma$ is the maximal $k$ such that there are indices $1 \leq i_1 < i_2 < \cdots < i_k \leq N$ satisfying $\sigma(i_1) < \sigma(i_2) < \cdots < \sigma(i_k)$. We denote this length by $\ell_N$. In the early 1960's, Ulam raised the question about the large-$N$ behavior of $\ell_N$. Based on computer simulations, he conjectured that $c = \lim_{N \to \infty} N^{-1/2} \mathbb{E}(\ell_N)$ exists in $(0, \infty)$. The verification of this statement and the identification of $c$ have become known as 'Ulam's problem'. A long list of researchers contributed to this problem, including Hammersley, Logan and Shepp, Vershik and Kerov, and Seppäläinen. By the end of the 1990's, it was known that the above limit exists with $c = 2$, and computer simulations suggested that[11]

$$\lim_{N \to \infty} \frac{\mathbb{E}(\ell_N) - 2\sqrt{N}}{N^{1/6}} \approx -1.7711. \tag{3.15}$$

A survey on the history of Ulam's problem may be found in [OR00] and [AD99].

There is a 'Poissonized' version of Ulam's problem, which is strongly related and provides a technical tool for the solution of Ulam's problem. Consider a homogeneous Poisson process on $(0, \infty)^2$ with parameter one, and let $L(\lambda)$ be the maximal number of points of this process which can be joined together by a polygon line that starts at $(0,0)$, ends at $(\sqrt{\lambda}, \sqrt{\lambda})$ and is always going in an up/right direction. Then it is easy to see[12] that the distribution of $L(\lambda)$ is equal to the distribution of $\ell_{N^*}$, where $N^*$ is a Poisson random variable with parameter $\lambda$. Via Tauberian theorems, asymptotics of the distribution of $L(\lambda)$ as $\lambda \to \infty$ stand in a one-to-one correspondence to the large-$N$ asymptotics of $\ell_N$.

There are exact formulas for the distributions both of $\ell_N$ and $L(\lambda)$, which have been proved by many authors using various methods (see [BDJ99]). Indeed, for any $n \in \mathbb{N}$, we have

$$\mathbb{P}(\ell_N \leq n) = \frac{2^{2N} N!}{(2N)!} \int_{[-\pi, \pi]^n} \Big(\sum_{j=1}^n \cos \theta_j\Big)^{2N} \prod_{1 \leq k < j \leq n} |e^{\mathrm{i}\, \theta_j} - e^{\mathrm{i}\, \theta_k}|^2 \frac{\mathrm{d}^n \theta}{(2\pi)^n n!},$$

$$\mathbb{P}(L(\lambda) \leq n) = e^{-\lambda} \int_{[-\pi, \pi]^n} \exp\Big\{2\sqrt{\lambda} \sum_{j=1}^n \cos \theta_j\Big\} \prod_{1 \leq k < j \leq n} |e^{\mathrm{i}\, \theta_j} - e^{\mathrm{i}\, \theta_k}|^2 \frac{\mathrm{d}^n \theta}{(2\pi)^n n!}. \tag{3.16}$$

In [BDJ99], sophisticated and deep methods are applied to the right hand side of (3.16), which have previously been established in [DZ93], [DZ95] and [DVZ97]:

---

[11] Recall Remark 2.18(iii).

[12] The main reason is the characteristic property of the Poisson process that, given that there are precisely $N$ Poisson points in the square, these points are conditionally independent and uniformly distributed.

the steepest-decent method for the computation of asymptotics of solutions to certain Riemann-Hilbert problems. As a result, a limit law for $\ell_N$ is proved, which shows again the universality of the Tracy-Widom distribution for GUE in (2.55):

**Theorem 3.14 (Limit law for $\ell_N$, [BDJ99]).** *Let $\ell_N$ be the length of the longest increasing subsequence of a random permutation, which is uniformly distributed over $\mathfrak{S}_N$. Then, as $N \to \infty$, the scaled random variable*

$$\chi_N = \frac{\ell_N - 2\sqrt{N}}{N^{1/6}} \tag{3.17}$$

*converges in distribution towards the Tracy-Widom distribution for GUE. Moreover, all moments of $\chi_N$ also converge towards the moments of this distribution. Both assertions are true also for $(L(\lambda) - 2\sqrt{\lambda})\lambda^{-\frac{1}{6}}$ as $\lambda \to \infty$.*

**Sketch of the proof.** We sketch some elements of the proof, partially also following [P03, Sect. 3.1]. We consider the Poissonized version and consider $L(\lambda^2)$ instead of $L(\lambda)$.

The starting point is an explicit expression for the probability of $\{L(\lambda^2) \leq N\}$ for any $N \in \mathbb{N}$ and any $\lambda > 0$ in terms of the Toeplitz determinant[13] $D_{N,\lambda} = \det T_N(e^{2\lambda \cos(\cdot)})$. More precisely, one has

$$\mathbb{P}(L(\lambda^2) \leq N) = e^{-\lambda^2} \det T_N(e^{2\lambda \cos(\cdot)}) = e^{-\lambda^2} D_{N,\lambda}, \tag{3.18}$$

a remarkable formula which has first been derived in [Ge90], based on the RSK-correspondence of Lemma 3.9. On $[0, 2\pi]$ we introduce the inner product

$$\langle p, q \rangle_\lambda = \int_0^{2\pi} p(e^{i\theta}) \overline{q(e^{i\theta})} e^{2\lambda \cos\theta} \frac{d\theta}{2\pi}. \tag{3.19}$$

Consider the sequence of orthogonal polynomials $(\pi_N^{(\lambda)})_{N \in \mathbb{N}_0}$ with respect to $\langle \cdot, \cdot \rangle_\lambda$ which is obtained via the Gram-Schmidt algorithm from the monomials $z^n$, $n \in \mathbb{N}_0$. We normalize $\pi_N^{(\lambda)}$ such that $\pi_N^{(\lambda)}(z) = z^N + \mathcal{O}(z^{N-1})$ and define $V_N^{(\lambda)} = \|\pi_N^{(\lambda)}\|^2$, such that we have

$$\langle \pi_N^{(\lambda)}, \pi_{N'}^{(\lambda)} \rangle_\lambda = \delta_{N,N'} V_N^{(\lambda)}, \qquad N, N' \in \mathbb{N}_0. \tag{3.20}$$

Classical results on orthogonal polynomials (see [Sz75] for some background) imply the identities

$$D_{N,\lambda} \equiv \det T_N(e^{2\lambda \cos(\cdot)}) = \det\bigl((\langle z^k, z^l \rangle_\lambda)_{k,l=0,\ldots,N-1}\bigr) = \prod_{k=0}^{N-1} V_k^{(\lambda)}$$
$$= (V_0^{(\lambda)})^N \prod_{k=0}^{N-1} \prod_{l=1}^{k} \bigl(1 - (\pi_l^{(\lambda)}(0))^2\bigr). \tag{3.21}$$

---

[13]We recall that the $(N \times N)$ Toeplitz matrix $T_N(f) = (\mu_{k-l})_{k,l=0,\ldots,N-1}$ with respect to the weight function $f$ on $[0, 2\pi]$ is defined by the Fourier coefficients $\mu_k = \int_0^{2\pi} e^{i k\theta} f(\theta) \frac{d\theta}{2\pi}$.

For our special choice of the weight function, $e^{2\lambda \cos \theta}$, one obtains a nonlinear recursion relation for the sequence $(\pi_N^{(\lambda)}(0))_{N \in \mathbb{N}_0}$, which are called the *discrete Painlevé II equations*. Indeed, the numbers $R_N^{(\lambda)} = (-1)^{N+1} \pi_N^{(\lambda)}(0)$ satisfy

$$R_{N+1}^{(\lambda)} - 2R_N^{(\lambda)} + R_{N-1}^{(\lambda)} = \frac{(\frac{N}{\lambda} - 2)R_N^{(\lambda)} + 2(R_N^{(\lambda)})^3}{1 - (R_N^{(\lambda)})^2}, \qquad N \in \mathbb{N}. \qquad (3.22)$$

Putting $N = \lfloor 2\lambda + \lambda^{\frac{1}{3}} s \rfloor$, multiplying (3.22) with $\lambda$ and letting $\lambda \to \infty$, we see that the function

$$\widetilde{q}(s) = - \lim_{\lambda \to \infty} \lambda^{\frac{1}{3}} R_{\lfloor 2\lambda + \lambda^{\frac{1}{3}} s \rfloor}^{(\lambda)}, \qquad s \in \mathbb{R}, \qquad (3.23)$$

should satisfy the (continuous) Painlevé II equation in (2.54). The initial value $R_0^{(\lambda)} = -1$, i.e., $\widetilde{q}(-2\lambda^{\frac{2}{3}}) \sim \lambda^{\frac{1}{3}}$, implies that we are dealing with that solution of (2.54) that is positive in $(-\infty, 0)$. Hence, $\widetilde{q}$ is identical to the solution $q$ of (2.54) with $q(x) \sim \mathrm{Ai}(x)$ as $x \to \infty$; recall the text below (2.54).

Note that (3.21) implies that $D_{N+1,\lambda} D_{N-1,\lambda} / D_{N,\lambda}^2 = 1 - (R_N^{(\lambda)})^2$. Using this in (3.18) we obtain, for $\lambda \to \infty$,

$$\left(\frac{\mathrm{d}}{\mathrm{d}s}\right)^2 \log \mathbb{P}\big((L(\lambda^2) - 2\lambda)\lambda^{-\frac{1}{3}} \leq s\big) \approx \lambda^{\frac{2}{3}} \Big(\log D_{N+1,\lambda} - 2\log D_{N,\lambda} + \log D_{N-1,\lambda}\Big)$$
$$= \lambda^{\frac{2}{3}} \log\Big(1 - \big(R_{\lfloor 2\lambda + \lambda^{\frac{1}{3}} s \rfloor}^{(\lambda)}\big)^2\Big) \approx -\Big(\lambda^{\frac{1}{3}} R_{\lfloor 2\lambda + \lambda^{\frac{1}{3}} s \rfloor}^{(\lambda)}\Big)^2 \approx -q(s)^2 = (\log F_2)''(s). \qquad (3.24)$$

Hence, we have finished the identification of the limiting distribution of $L(\lambda^2)$.

The technically hardest works of the proof are the proofs of the convergence in (3.23) and of the convergence of the moments, which require an adaptation of the Deift-Zhou steepest descent method for an associated Riemann-Hilbert problem. □

### 3.6 The poly nuclear growth model

Consider the boundary of a one-dimensional substrate, which is formed by the graph of a piecewise constant function with unit steps. At each time $t \geq 0$, the separation line between the substrate and its complement is given as the graph of the function $h(\cdot, t) \colon \mathbb{R} \to \mathbb{R}$. Occasionally, there occur random nuclear events in states $x^*$ at times $t^*$, and the process of the pairs $(x^*, t^*)$ forms a Poisson point process in the space-time half plane $\mathbb{R} \times [0, \infty)$ with intensity equal to two. Such an event creates an island of height one with zero width, i.e., $h$ has a jump of size one in $x^*$ at time $t^*$. Every island grows laterally (deterministically) in both directions with velocity one, but keeps its height, i.e., for small $\varepsilon > 0$ the curve $h(\cdot, t^* + \varepsilon)$ has the height $h(x^*, t^*)$ in the $\varepsilon$-neighborhood of $x^*$ and stays on the same level as before $t^*$ outside this neighborhood:

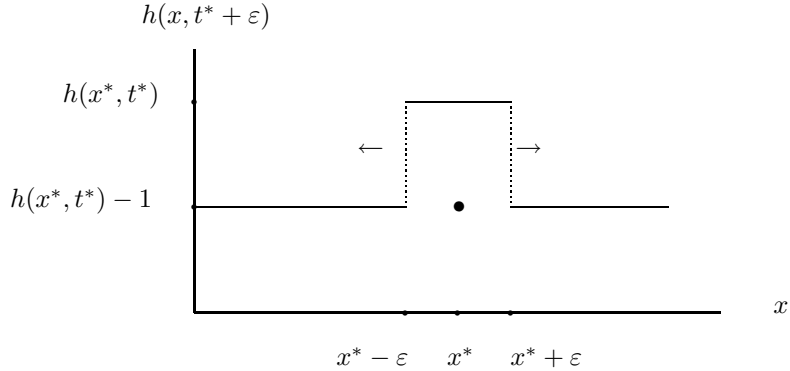

The bullet marks the nucleation event, and the two arrows indicate the lateral growth of velocity one in the two directions. We call the graph of $h(\cdot, t^* + \varepsilon)$ in the $\varepsilon$-neighborhood of $x^*$ a *growing island*. If two growing islands at the same level collide, then they merge together and form a common growing island. The nucleation events occur only on top of a growing island, and they occur with constant density equal to two.

This is a (rather simple) model for *poly nuclear growth* (PNG) in $1+1$ dimension. Among various initial conditions that one could impose, we shall consider only two: the *flat* case, where $h(x, 0) = 0$ for any $x \in \mathbb{R}$, and the *droplet* case, where $h(x, 0) = -\infty$ for $x \neq 0$ and $h(0, 0) = 0$. The droplet case may also be defined with the initial condition $h(\cdot, 0) = 0$ by requiring that nucleation events at time $t$ may happen only in $[-t, t]$.

Let us first consider the droplet case. A beautiful observation [PS00] is the fact that the PNG model stands in a one-to-one relation to the Poissonized problem of the longest increasing subsequence in a rectangle. Using this correspondence, one arrives at the following limit assertion.

**Theorem 3.15 (Limit law for the PNG droplet, [PS00]).** *Let $h(x, t)$ be the height of the PNG droplet at time $t$ over the site $x$, and let $c \in [-1, 1]$. Then*

$$\lim_{N \to \infty} \mathbb{P}\Big( \frac{h(ct, t) - 2t\sqrt{1-c^2}}{(1-c^2)^{\frac{1}{6}} t^{\frac{1}{3}}} \leq s \Big) = F_2(s), \qquad s \in \mathbb{R}, \qquad (3.25)$$

*where $F_2$ is the GUE Tracy-Widom distribution function, see* (2.55).

**Idea of proof.** We consider the space-time half plane $\mathbb{R} \times [0, \infty)$. For any space-time point $(x, t)$, we call the quarter plane with lower corner at $(x, t)$ and having the two lines through $(x, t)$ with slope 1 and $-1$ as boundaries the $(x, t)$-quarter plane. Recall that nucleation events occur in the $(0, 0)$-quarterplane only, which is the region $\{(x, t): |x| \leq t\}$.

First note that every nucleation event at some space-time point $(x^*, t^*)$ influences the height of the curve $h$ only within the $(x^*, t^*)$-quarter plane. Second, note that any nucleation event $(y^*, s^*)$ within the $(x^*, t^*)$-quarter plane contributes an additional lifting by level one (to the lift created by the nucleation event $(x^*, t^*))$ for any space-time point in the intersection of the two quarter planes of the nucleation events, since the growing island created by $(y^*, s^*)$ will be on top of the growing island created by $(x^*, t^*)$. However, if $(y^*, s^*)$ occurs outside the $(x^*, t^*)$-quarter plane, their influences are merged to a lift just by one step since their growing islands are merged to one growing island.

Now fix a space-time point $(x, t)$ in the $(0, 0)$-quarter plane. In the space-time plane, consider the rectangle $R$ having two opposite corners at the origin and at the point $(x, t)$ and having sides of slopes 1 and $-1$ only. Condition on a fixed number $N$ of nucleation events $(x_1^*, t_1^*), \ldots, (x_N^*, t_N^*)$ in the rectangle $R$.

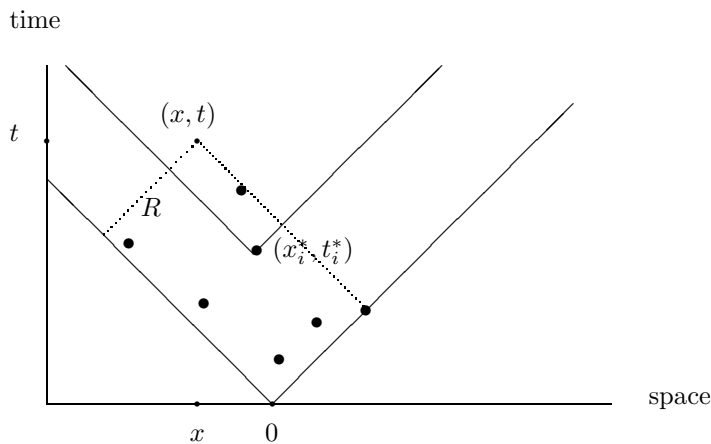

Rotate the rectangle by 45 degrees. The preceding observations imply that only those nucleation events contribute to the height $h(x, t)$ which can be joined together by a polygon line consisting of straight up/right lines, leading from the corner of the rectangle $R$ at the origin to the corner at $(x, t)$. The maximal number of nucleation events along such path is equal to the height $h(x, t)$. Hence, the length of the longest increasing subsequence in a unit square with Poisson intensity $\lambda = \sqrt{t^2 - x^2}$ has the same distribution as the height $h(x, t)$. Using Theorem 3.14, one concludes the assertion. □

In particular, the fluctuation exponent $1/3$ is rigorously proved for this growth model. Such a result has not yet been achieved for any other growth model of this type. However, this fluctuation behavior is conjectured for a large class of $(1 + 1)$-dimensional growth processes, provided the spatial correlations are not too weak.

The flat initial condition, $h(\cdot, 0) = 0$, interestingly leads to the GOE Tracy-Widom distribution instead of the GUE one:

**Theorem 3.16 (Limit law for the flat PNG model, [PS00]).** *Let $h(x, t)$ be the height of the flat PNG model at time $t$ over the site $x$. Then,*

$$\lim_{t \to \infty} \mathbb{P}\Big(\frac{h(0, t) - 2t}{t^{\frac{1}{3}}} \leq 2^{-2/3} s\Big) = F_1(s), \qquad s \in \mathbb{R}, \qquad (3.26)$$

*where $F_1$ is the GOE Tracy-Widom distribution function, see* (2.58).

The above explanation for the droplet case has to be adapted to the flat case by replacing the *rectangle* with corners at the origin and $(x,t)$ by the *triangle* with base on the axis $t = 0$, corner at $(x,t)$ and side slopes 1 and $-1$. See [Fe04a] for more detailed results on the flat PNG model.

For other initial conditions (among which some lead to the GSE Tracy-Widom distribution, $F_4$), see [P03, Sect. 3]. We recall that a discrete-space version of the PNG model is analyzed in [Jo02a, Sect. 2.4]; see also Remark 3.11. A recent survey on the PNG droplet and its relation to random matrices and further random processes, like directed polymers, the longest increasing subsequence problem and Young tableaux, appears in [FP05].

### 3.7 The multi-layer PNG droplet and the Airy process

The PNG droplet has been analysed also as a process. Interestingly, the limiting distribution of the height process in the correct scaling bears a close relationship to Dyson's Brownian motions (see Theorem 4.1), which is best seen when additional layers of substrate separation lines are introduced. The so-called *multi-layer PNG droplet* (sometimes also called the *discrete PNG model*) is defined as follows. We write $h_0$ instead of $h$ and add an infinite sequence of separation lines $h_\ell(x,t)$ with $\ell \in -\mathbb{N}$ with initial condition $h_\ell(x,0) = \ell$. Nucleation events only occur to the zeroth line $h_0$, and they occur at time $t$ in the interval $[-t,t]$ only (i.e., we consider the droplet case). Every merging event in the $\ell$-th line (i.e., every event of an amalgamation of two neighboring growing islands at the same height) creates a nucleation event in the $(\ell - 1)$-st line at the same site. Apart from this rule, every island on any level grows deterministically with unit speed into the two lateral directions as before.

Hence, randomness is induced only at the zeroth line, and all the other lines are deterministic functions of $h_0$. Observe that the strict ordering $h_\ell(x,t) > h_{\ell-1}(x,t)$ for any $x,t,\ell$ is preserved. Hence, the lines form a family of non-colliding step functions with unit steps. For any $\ell \in -\mathbb{N}_0$ and at any time $t > 0$, the $\ell$-th line $h_\ell(\cdot,t)$ is constant equal to $\ell$ far away from the origin. Only a finite (random) number of them have received any influence coming from the nucleation events, and only within a finite (random) space-time window.[14]

An interesting observation [PS02a] is that, in the long-time limit, the multi-layer PNG droplet process approaches the large-$N$ limit of Dyson's Brownian motions (see Section 4.1 below) in the appropriate scaling.[15] More precisely, let $\lambda^{(N)}(t) = (\lambda_1^{(N)}(t), \ldots, \lambda_N^{(N)}(t)) \in W_N$ be Dyson's Brownian motion at time $t$ as in Theorem 4.1. Then the *Airy process* may be introduced as the scaled limiting distribution of the largest particle, more precisely,

$$\left(\sqrt{2}N^{\frac{1}{6}}\left(\lambda_N^{(N)}(yN^{-\frac{1}{3}}) - \sqrt{2N}\right)\right)_{y \in \mathbb{R}} \implies \left(\mathfrak{Ai}(y)\right)_{y \in \mathbb{R}}. \qquad (3.27)$$

---

[14]Computer simulations show that the space-time region in which the lines are not constant asymptotically forms a droplet that approaches a circle. This region stays strictly inside the circle, which is due to the negativity of the expectation of the GUE Tracy-Widom distribution, recall Remark 2.18(iii).

[15]To be more precise, in contrast to Section 4.1, here Dyson's Brownian motions are not based on Brownian motions, but on Ornstein-Ulenbeck processes, which are Brownian motions with drift to the origin and hence stationary.

Convergence has been established in the sense of finite-dimensional distributions in [PS02a] and in process sense in [Jo03]. For any $y > 0$, the random variable $\mathfrak{A}i(y)$ has the GUE Tracy-Widom distribution $F_2$ in (2.55), and the family of these random variables forms an interesting stochastic process. The *Airy process* $(\mathfrak{A}i(y))_{y \in \mathbb{R}}$ is a stationary, continuous non-Markovian stochastic process which may be defined via its finite dimensional distributions, using a determinant formula involving the Airy kernel $K_{\text{Ai}}$ in (2.52) [PS02a], see also [P03, Sect. 5].

In [PS02a] it turns out that, in the appropriate scaling, the joint distribution of all the lines $h_\ell$ of the multilayer PNG droplet approaches the Airy process. We state the consequence of this statement for the first line as follows.

**Theorem 3.17 (Process convergence of the PNG droplet, [PS02a]).**
*Let $h(x, t)$ be the height of the PNG droplet at time t over the site x. Then, in the sense of finite-dimensional distributions,*

$$\lim_{t \to \infty} \frac{h(t^{\frac{2}{3}} y, t) - 2t}{t^{\frac{1}{3}}} = \mathfrak{A}i(y) - y^2, \qquad y \in \mathbb{R}, \tag{3.28}$$

*where $(\mathfrak{A}i(y))_{y \in \mathbb{R}}$ is the Airy process.*

Some progress on the process version of the *flat* PNG model has been made in [Fe04a]. Discrete versions of the PNG model have been analysed in [IS04a], [IS04b].

Another interesting process that converges (after proper rescaling) in distribution towards the Airy process is the north polar region of the *Aztec diamond* [Jo05].

### 3.8 The Plancherel measure

The Plancherel measure is a distribution on the set of Young tableaux which exhibits an asymptotic behavior that is remarkably similar to that of the spectrum of Gaussian matrix ensembles. Most interestingly, this measure may be studied for any value of the parameter $\beta$, which is restricted to the values 1, 2 and 4 in the matrix cases.

A *Young diagram*, or equivalently a partition $\lambda = (\lambda_1, \lambda_2, \dots)$ of $\{1, \dots, N\}$ is an array of $N$ boxes, such that $\lambda_1$ of them are in the first row, $\lambda_2$ of them in the second and so on. Here $\lambda$ is an integer-valued partition such that $\lambda_1 \geq \lambda_2 \geq \dots$, and $\sum_i \lambda_i = N$. We think of the rows as being arranged on top of each other. A *standard Young tableau* is a Young diagram together with a filling of the boxes with the numbers $1, \dots, N$ such that the numbers are strictly increasing along the rows and along the columns.[16] The vector $\lambda$ is called the *shape* of the tableau. For every $\lambda$, we denote by $d_\lambda$ the number of Young tableaux of shape $\lambda$. For every $\beta > 0$, we define the *Plancherel measure* as the distribution on the set $Y_N$ of partitions of $\{1, \dots, N\}$, which is given by

$$\text{Pl}_N^{(\beta)}(\lambda) = \frac{d_\lambda^\beta}{\sum_{\mu \in Y_N} d_\mu^\beta}, \qquad \lambda \in Y_N. \tag{3.29}$$

---

[16]Compare to the definition of a semistandard Young tableau prior to Lemma 3.9, where more numbers may appear, and their order is just nondecreasing.

We can conceive $\lambda_k^{(N)}$, the length of the $k$-th row, as an $\mathbb{N}_0$-valued random variable under the probability measure $\mathrm{Pl}_N^{(\beta)}$ on $Y_N$.

The case $\beta = 2$ has been studied a lot. Basically, it was shown that the limiting statistics of the sequence $\lambda_1^{(N)}, \lambda_2^{(N)}, \ldots$, in an appropriate scaling, is the same as the one for the eigenvalues of an $(N \times N)$ GUE-matrix. We mention just a few important results. As a by-product of their study of the longest increasing subsequence of a random permutation, in [BDJ99] the limit theorem

$$\lim_{N\to\infty} \mathrm{Pl}_N^{(2)}\Big(\frac{\lambda_1^{(N)} - 2\sqrt{N}}{N^{\frac{1}{6}}} \leq s\Big) = F_2(s), \qquad s \in \mathbb{R}, \qquad (3.30)$$

is shown, where $F_2$ is the Tracy-Widom GUE distribution function. The conjecture of [BDJ99] that for every $k \in \mathbb{N}$ the scaled limiting distribution of $\lambda_k^{(N)}$ is identical to the one of the $k$-th largest eigenvalue of a GUE-matrix was independently proved in [BDJ00] for $k = 2$, and for general $k$ in [Jo01b] and [BOO00]. The convergence of the moments of the scaled row lengths was also proved in [BDJ99], [BDJ00] and [Jo01b], respectively. The bulk-scaling limit was also proved in [BOO00]. The case $\beta = 1$ (which is analogous to the GOE case instead the GUE case) has been studied in [BE01].

## 4. Non-colliding random processes

In this section we systematically discuss conditional multi-dimensional random processes given that the components never collide with each other. These processes are sometimes called *vicious walkers*, *non-colliding processes* or *nonintersecting paths* in the literature. The earliest hint at a close connection between non-colliding random processes and orthogonal polynomial ensembles was found in [Dy62b], where a natural process version of the Gaussian Unitary Ensemble was considered. It turned out there that the mutual repellence in (2.3) receives a natural interpretation in terms of Brownian motions conditioned on never colliding with each other. This theme apparently was not taken up in the literature up to the beginning of the nineties, when people working in stochastic analysis turned to this subject. Since the discovery of close connections also with random growth models at the end of the nineties, non-colliding processes became an active research area.

### 4.1 Dyson's Brownian motions

A glance at the Hermite ensemble in (2.3) shows that there is a mutually repelling force between the eigenvalues: the density vanishes if any two of the $N$ arguments approach each other. It does not seem easy to derive an intuitive reason for this repellence from random matrix considerations, but if the matrix $M$ is embedded in a natural process of random Hermitian matrices, then the process of eigenvalues admits a nice identification that makes the repellence natural.

**Theorem 4.1 (Dyson's Brownian motions, [Dy62b]).** *For any $i \in \{1,\ldots,N\}$ resp. $i < j$, let $(M_{i,i}(t))_{t\geq 0}$ and $(M_{i,j}^{(\mathrm{R})}(t))_{t\geq 0}$ and $(M_{i,j}^{(\mathrm{I})}(t))_{t\geq 0}$ be* *independent real standard Brownian motions, starting at zero, such that the Hermitian random matrix $M(t) = (M_{i,j}(t))_{i,j=1,\ldots,N}$ with $M_{i,j}(t) = M_{i,j}^{(\mathrm{R})}(t) + \mathrm{i}\, M_{i,j}^{(\mathrm{I})}(t)$ has the distribution of GUE at time $t = 1$. Then the process $(\lambda(t))_{t\geq 0}$ of eigenvalue vectors $\lambda(t)$ of $M(t)$ is a conditional Brownian motion on $\mathbb{R}^N$, starting at zero, given that the particles never collide with each other, i.e., conditional on the event $\{\lambda_1(t) < \lambda_2(t) < \cdots < \lambda_N(t) \text{ for all } t > 0\}$.*

This theorem has to be explained in some detail.

**Remark 4.2.**  (i) It is remarkable that, in particular, the process of eigenvalue vectors is Markov. This is not true for, say, the process of the largest eigenvalues, $(\lambda_N(t))_{t\geq 0}$.

(ii) The original proof in [Dy62b] makes nowadays an old-fashioned impression. See [Br91] for a modern stochastic analysis treatment of an analogous matrix-valued process for Wishart-matrices in the real-valued setting. In this setting, the process of eigenvalues also turns out to be Markov, but does not admit a conditional interpretation. The latter is also true in the analogous GOE setting.

(iii) The event of never colliding, $\{\lambda_1(t) < \lambda_2(t) < \cdots < \lambda_N(t) \text{ for all } t > 0\}$, has zero probability for $N$ independent Brownian motions. Hence, the definition of the conditioned process needs some care. First observe that the non-colliding event is the event $\{\lambda(t) \in W_N \text{ for all } t > 0\}$, where $W_N = \{x \in \mathbb{R}^N \colon x_1 < x_2 < \cdots < x_N\}$ is the Weyl chamber. Probabilists like to write this event as $\{T = \infty\}$, where $T = \inf\{t > 0 \colon \lambda(t) \in W_N^{\mathrm{c}}\}$ is the exit time from $W_N$, the first time of a collision of any two of the particles. One way to construct the conditional process is to condition on the event $\{T > t\}$ and prove that there is a limiting process as $t \to \infty$. Another one is to consider the Doob-$h$ transform of the vector of $N$ independent standard Brownian motions with some suitable function $h \colon W_N \to (0, \infty)$ that vanishes on the boundary of $W_N$ and is harmonic for the generator of the $N$-dimensional Brownian motion in $W_N$. Remarkably, it turns out that $h = \Delta_N$, the Vandermonde determinant, satisfies all these properties, and that the $h$-transform with this function $h$ is identical with the outcome of the first construction. See Section 4.2 below for a general treatment of this issue.

(iv) The Markov process $(\lambda(t))_{t\geq 0}$ has the invariant measure $x \mapsto \Delta_N(x)^2 \,\mathrm{d}x$, which cannot be normalized.

(v) Also in the real and the symplectic version, the eigenvalue process, $(\lambda(t))_{t\geq 0}$, turns out to be a diffusion. An elementary application of Ito's formula shows that $(\lambda(t))_{t\geq 0}$ satisfies the stochastic differential equation (see [Br91] for related formulas)

$$\mathrm{d}\lambda_i = \mathrm{d}B_i + \beta \sum_{\substack{j=1 \\ j\neq i}}^{N} \frac{1}{\lambda_i - \lambda_j}\,\mathrm{d}t, \qquad i = 1,\ldots,N, \tag{4.1}$$

where $B_1, \ldots, B_N$ are independent Brownian motions, and $\beta \in \{1, 2, 4\}$ is the parameter as in (2.5). The generator of the process $(\lambda(t))_{t\geq 0}$ is

$$\widehat{G}f(x) = \frac{1}{2}\sum_{i=1}^{N}\partial_i^2 f(x) + \beta \sum_{i=1}^{N}\Big[\sum_{\substack{j=1\\j\neq i}}^{N}\frac{1}{x_i - x_j}\Big]\partial_i f(x), \qquad (4.2)$$

The generators in the GOE and the GSE setting have a factor different from 2 before the drift term. Apparently this circumstance makes it impossible to conceive the processes as Doob transforms of $N$ independent processes. ◇

## 4.2 Harmonicity of the Vandermonde determinant

Now we consider more general multi-dimensional random processes and their conditional version given that no collision of the particles occurs. As was pointed out in Remark 4.2(iii), the construction needs some care, since the conditioning is on a set of probability zero. It turns out that the rigorous definition may be given for many processes in terms of a Doob $h$-transform with $h = \Delta_N$, the Vandermonde determinant in (1.3). Even more striking, the marginal distribution of the transformed process turns out, for a couple of explicit examples, to be given by well-known orthogonal polynomial ensembles, like the Hermite ensemble in (2.3) for the case of conditional Brownian motions.

*4.2.1. The continuous case.*

Let us turn first to the time-continuous case with continuous paths, more precisely, to diffusions. We fix $N \in \mathbb{N}$ and an interval $I$ and let $X = (X(t))_{t\geq 0}$ be a stochastic process on $I^N$. Assume that $X_1, \ldots, X_N$ are $N$ independent and identically distributed diffusions $X_i = (X_i(t))_{t\geq 0}$ on $I$. Under the measure $\mathbb{P}_x$ they start at $X_i(0) = x_i \in I$, where $x = (x_1, \ldots, x_N)$. By $p_t(x,y)$ we denote the transition density function of any of the diffusions $X_i$, i.e.,

$$\mathbb{P}_x(X(t) \in dy) = \prod_{i=1}^{N}\big[p_t(x_i, y_i)\, dy_i\big], \qquad x, y \in I^N. \qquad (4.3)$$

Recall the Weyl chamber and its exit time,

$$W_N = \{x \in \mathbb{R}^N : x_1 < \cdots < x_N\} \qquad \text{and} \qquad T = \inf\{t > 0 : X(t) \notin W_N\}. \qquad (4.4)$$

In words: $T$ is the first time of a collision of any two of the $N$ components of the process. Recall the Vandermonde determinant $\Delta_N(x) = \prod_{1\leq i<j\leq N}(x_j - x_i)$. In order to be able to construct a Doob-$h$ transform of the process with $h = \Delta_N$ on $W_N$, the basic requirements are: (1) $\Delta_N$ is positive on $W_N$, (2) $\Delta_N$ is harmonic with respect to the generator $G$ of the process $X$, i.e., $G\Delta_N = 0$, and (3) $\Delta_N(X(t))$ is integrable for any $t > 0$.

Clearly, the first prerequisite is satisfied. Furthermore, it turns out that $\Delta_N$ is harmonic for a quite large class of processes:

**Lemma 4.3 (Harmonicity of $\Delta_N$, continuous case [KO01]).** *We have $G\Delta_N = 0$ (i.e., $\Delta_N$ is harmonic with respect to $G$) if there are $a, b, c \in \mathbb{R}$ such that*

$$G = \sum_{i=1}^{N}\Big[(ax_i+b)\partial_i^2+c\partial_i\Big] \quad or \quad = \sum_{i=1}^{N}\Big[(x_i^2+ax_i+b)\partial_i^2+\big(\tfrac{2}{3}(N-2)x_i+c\big)\partial_i\Big]. \tag{4.5}$$

The proof consists of an elementary calculation. Lemma 4.3 in particular covers the cases of Brownian motion, squared Bessel processes (squared norms of Brownian motions) and generalized Ornstein-Uhlenbeck processes driven by Brownian motion. For general diffusions, existence and identification of positive harmonic functions for the restriction of the generator to the Weyl chamber are open.

As a consequence of Lemma 4.3, we can introduce the Doob $h$-transform of $X$ with $h = \Delta_N$. This is a diffusion on $W_N \cap I^N$, which we also denote $X$. Its transition probability function is given by

$$\widehat{\mathbb{P}}_x(X(t) \in \mathrm{d}y) = \mathbb{P}_x(T > t; X(t) \in \mathrm{d}y)\frac{\Delta_N(y)}{\Delta_N(x)}, \qquad x, y \in W_N \cap I^N, t > 0. \tag{4.6}$$

The transformed process is often called the conditional process $X$, given that there is no collision of the components. In order to justify this name, one must show that

$$\lim_{t\to\infty}\mathbb{P}_x(X(s) \in \mathrm{d}y \mid T > t) = \widehat{\mathbb{P}}_x(X(s) \in \mathrm{d}y), \qquad \text{for any } x, y \in W_N, s > 0. \tag{4.7}$$

This may be proven in many examples with the help of the Markov property at time $s$ and an asymptotic formula for $\mathbb{P}_z(T > t)$ as $t \to \infty$, see Remark 4.10(ii). In Section 4.3 we provide two tools. In Section 4.4, we list a couple of examples of $\Delta_N$-transformed diffusions, whose marginal distribution is an orthogonal polynomial ensemble.

*4.2.2. The discrete case.*

There is also a discrete version of Lemma 4.3. Recall that a vector $v$ on a discrete set $I$ is called a *positive regular function* for a matrix $Q$ with index set $I \times I$ if all the components of $v$ are positive and $Qv = v$ holds.

**Lemma 4.4 (Regularity of $\Delta_N$, discrete case [KOR02]).** *Let $(X(n))_{n\in\mathbb{N}}$ be a random walk on $\mathbb{R}^N$ such that the step distribution is exchangeable and the $N$-th moment of the steps is finite.*

(i) *Then $\Delta_N$ is harmonic for the walk, i.e., $\mathbb{E}_x[\Delta_N(X(1))] = \Delta_N(x)$ for any $x \in \mathbb{R}^N$, and the process $\Delta_N(X(n))_{n\in\mathbb{N}_0}$ is a martingale with respect to the natural filtration of $(X(n))_{n\in\mathbb{N}}$.*

(ii) If $(X(n))_n$ takes values in $\mathbb{Z}^N$ only and no step from $W_N$ to $\overline{W}_N^c$ has positive probability, then the restriction of $\Delta_N$ to $W_N \cap \mathbb{Z}^N$ is a positive regular function for the restriction $P_{W_N} = (p(x,y))_{x,y \in W_N \cap \mathbb{Z}^N}$ of the transition matrix $P = (p(x,y))_{x,y \in \mathbb{Z}^N}$, i.e.,

$$\sum_{y \in \mathbb{Z}^N \cap W_N} p(x,y) \Delta_N(y) = \Delta_N(x), \qquad \text{for any } x \in \mathbb{Z}^N \cap W_N. \quad (4.8)$$

The condition in Lemma 4.4(ii) is a severe restriction. It in particular applies to nearest-neighbor walks on $\mathbb{Z}^N$ with independent components, and to the multinomial walk, where at each discrete time unit one randomly chosen component makes a unit step, see Section 4.4. Further examples comprise birth and death processes and the Yule process [Do05, Ch. 6].

Under the assumptions of Lemma 4.4, one can again define the $h$-transform of the Markov chain $X$ by using the transition matrix $\widehat{P} = (\widehat{p}(x,y))_{x,y \in W_N \cap \mathbb{Z}^N}$ with

$$\widehat{p}(x,y) = p(x,y) \frac{\Delta_N(y)}{\Delta_N(x)}, \qquad x, y \in W_N \cap \mathbb{Z}^N.$$

**Remark 4.5.** *Arbitrary* random walks with i.i.d. components are considered in [EK05+]. Under the sole assumption of finiteness of sufficiently high moments of the steps, it turns out there that the function

$$V(x) = \Delta_N(x) - \mathbb{E}_x[\Delta_N(X(\tau))], \qquad x \in W_N,$$

where $\tau = \inf\{n \in \mathbb{N} \colon X(\tau) \notin W_N\}$ is the exit time from $W_N$, is a positive regular function for the restriction of the walk to $W_N$. (Note that $V$ coincides with $\Delta_N$ in the special cases of Lemma 4.4(ii).) Since the steps are now arbitrarily large, the term 'non-colliding' should be replaced by 'ordered'. Furthermore, an ordered version of the walk is constructed in terms of a Doob $h$-transform with $h = V$, and some asymptotic statements are derived, in particular an invariance principle towards Dyson's Brownian motions. ◇

### 4.3 Some tools

We present two technical tools that prove useful in the determination of probabilities of non-collision events.

*4.3.1. The Karlin-McGregor formula*

An important tool for calculating non-colliding probabilities is the *Karlin-McGregor formula*, which expresses the marginal distribution of the non-colliding process in terms of a certain determinant.

**Lemma 4.6 (Karlin-McGregor formula, [KM59]).** *Let $(X(t))_{t \geq 0}$ be a diffusion on $\mathbb{R}^N$ that satisfies the strong Markov property. Then, for any $x, y \in W_N$ and any $t > 0$,*

$$\frac{\mathbb{P}_x(T > t, X(t) \in \mathrm{d}y)}{\mathrm{d}y} = \det\big[(p_t(x_i, y_j))_{i,j=1,\ldots,N}\big], \quad (4.9)$$

where $p_t(x,y)$ is the transition probability function of the diffusion, see (4.3).

**Proof.** By $\mathfrak{S}_N$ we denote the set of permutations of $1,\ldots,N$, and $\text{sign}(\sigma)$ denotes the signum of a permutation $\sigma$. We write $y_\sigma = (y_{\sigma(1)},\ldots,y_{\sigma(N)})$. We have

$$\frac{\mathbb{P}_x(T > t, X(t) \in dy)}{dy} - \det\left[(p_t(x_i,y_j))_{i,j=1,\ldots,N}\right]$$
$$= \sum_{\sigma \in \mathfrak{S}_N} \text{sign}(\sigma)\left[\frac{\mathbb{P}_x(T > t, X(t) \in dy_\sigma)}{dy} - \frac{\mathbb{P}_x(X(t) \in dy_\sigma)}{dy}\right]$$
$$= -\sum_{\sigma \in \mathfrak{S}_N} \text{sign}(\sigma)\frac{\mathbb{P}_x(T \leq t, X(t) \in dy_\sigma)}{dy}, \tag{4.10}$$

since all the summands $\mathbb{P}_x(T > t, X(t) \in dy_\sigma)/dy$ are equal to zero, with the exception of the one for the identical permutation.

At time $T$, the $i$-th and the $j$-th coordinate of the process coincide for some $i < j$, which we may choose minimal. Reflect the path $(X(s))_{s \in [T,t]}$ in the $(i,j)$-plane, i.e., map this path onto the path $(X_\lambda(s))_{s \in [T,t]}$, where $\lambda \in \mathfrak{S}_N$ is the transposition that interchanges $i$ and $j$. This map is measure-preserving, and the endpoint of the outcome is at $y_{\sigma \circ \lambda}$ if $X(t) = y_\sigma$. Summing on all $i < j$ (i.e., on all transpositions $\lambda$), substituting $\sigma \circ \lambda$ and noting that its signum is the negative signum of $\sigma$, we see that the right hand side of (4.10) is equal to its negative value, i.e., it is equal to zero. The proof is finished. $\square$

**Remark 4.7.**   (i) The main properties of the process that make this proof possible are the strong Markov property and the continuity of the paths. No assumption on spatial dependence of the transition probability function is needed.

 (ii) For discrete-time processes on $\mathbb{Z}$ there is an analogous variant of Lemma 4.6, but a kind of continuity assumption has to be imposed: The steps must be $-1$, $0$ or $1$ only, i.e., it must be a nearest-neigbor walk. This ensures that the path steps on the boundary of $W_N$ when leaving $W_N$, and hence the reflection procedure can be applied. $\diamond$

*4.3.2. The Schur polynomials*

Another useful tool when dealing with certain determinants is the *Schur polynomial*,

$$\text{Schur}_z(x) = \frac{\det\left[(x_i^{z_j})_{i,j=1,\ldots,N}\right]}{\Delta_N(x)}, \qquad z \in W_N, x \in \mathbb{R}^N. \tag{4.11}$$

It turns out that $\text{Schur}_z$ is a multipolynomial in $x_1,\ldots,x_N$, and it is homogeneous of degree $z_1 + \cdots + z_N - \frac{N}{2}(N-1)$. Its coefficients are nonnegative integers and may be defined in a combinatorial way. It has the properties $\text{Schur}_z(1,\ldots,1) = \Delta_N(z)/\Delta_N(x^*)$ (where we recall that $x^* = (0,1,2,\ldots,N-1)$), $\text{Schur}_{x^*}(x) = 1$ for any $x \in \mathbb{R}^N$, and $\text{Schur}_z(0,\ldots,0) = 0$ for any $z \in W_N \setminus \{x^*\}$.

A combination of the Karlin-McGregor formula and the Schur polynomials identifies the asymptotics of the non-collision probability and the limiting joint distribution of $N$ standard Brownian motions before the first collision:

**Lemma 4.8.** *Let $(X(t))_{t \geq 0}$ be a standard Brownian motion, starting at $x \in W_N$. Then, as $t \to \infty$, the limiting distribution of $t^{-\frac{1}{2}} X(t)$ given that $T > t$ has the density $y \mapsto \frac{1}{Z} \varphi(y) \Delta_N(y)$ on $W_N$, where $\varphi$ is the standard Gaussian density, and $Z$ the normalization constant. Furthermore, $\mathbb{P}_x(T > t) = \Delta_N(x) t^{-\frac{N}{4}(N-1)} (C + o(1))$ as $t \to \infty$ for some $C > 0$.*

Note that the limiting distribution is of the form (1.1) with $\Delta_N^2$ replaced by $\Delta_N$, i.e., with $\beta = 1$.

**Sketch of proof.** Lemma 4.6 yields

$$\frac{\mathbb{P}_x(T > t, t^{-\frac{1}{2}} X(t) \in dy)}{dy} = \det\left[\left((2\pi)^{-\frac{N}{2}} e^{-(x_i - y_j\sqrt{t})^2/(2t)}\right)_{i,j=1,\ldots,N}\right]$$
$$= (2\pi)^{-\frac{N}{2}} e^{-\|x\|_2^2/(2t)} e^{-\|y\|_2^2/2} \det\left[\left(e^{x_i y_j/\sqrt{t}}\right)_{i,j=1,\ldots,N}\right]$$
$$= \frac{e^{-\|y\|_2^2/2}}{(2\pi)^{\frac{N}{2}}} e^{-\|x\|_2^2/(2t)} \Delta_N(z) \text{Schur}_y(z), \quad (4.12)$$

where we put $z_i = e^{x_i/\sqrt{t}}$. Now we consider the limit as $t \to \infty$. The second term is $(1 + o(1))$, and the continuity of $\text{Schur}_y$ implies that the last term converges to $\Delta_N(y)/\Delta_N(x^*)$. Using the approximation $e^{x_i/\sqrt{t}} - 1 \sim x_i/\sqrt{t}$, we see that $\Delta_N(z) \sim t^{-\frac{N}{4}(N-1)} \Delta_N(x)$. Hence, the right hand side of (4.12) is equal to $\varphi(y) \Delta_N(y) \Delta_N(x) t^{-\frac{N}{4}(N-1)} (1/\Delta_N(x^*) + o(1))$. Integrating on $y \in W_N$, we obtain the last statement of the lemma. Dividing the left hand side of (4.12) by $\mathbb{P}_x(T > t)$ and using the above asymptotics, we obtain the first one. □

### 4.4 Marginal distributions and ensembles

We apply now the technical tools of Section 4.3 to identify the marginal distribution of some particular $\Delta_N$-transformed processes as certain orthogonal polynomial ensembles.

*4.4.1. The continuous case.*

**Lemma 4.9 (Marginal distribution for $\Delta_N$-transformed diffusions, [KO01]).** *Assume that $I$ is an interval and $X$ is a diffusion on $I^N$ such that the Vandermonde determinant $\Delta_N$ is harmonic for its generator and $\Delta_N(X(t))$ is integrable for any $t > 0$. Assume that there is a Taylor expansion*

$$\frac{p_t(x,y)}{p_t(0,y)} = f_t(x) \sum_{m=0}^{\infty} (xy)^m a_m(t), \qquad t \geq 0, y \in I,$$

for $x$ in a neighborhood of zero, where $a_m(t) > 0$ and $f_t(x) > 0$ satisfy $\lim_{t\to\infty} a_{m+1}(t)/a_m(t) = 0$ and $f_t(0) = 1 = \lim_{t\to\infty} f_t(x)$. Then, for any $t > 0$ and some suitable $C_t > 0$,

$$\lim_{\substack{x \to 0 \\ x \in W_N}} \widehat{\mathbb{P}}_x(X(t) \in \mathrm{d}y) = C_t \Delta_N(y)^2 \mathbb{P}_0(X(t) \in \mathrm{d}y), \qquad y \in W_N. \quad (4.13)$$

Furthermore, for any $x \in W_N$,

$$\mathbb{P}_x(T > t) \sim C_t \Delta_N(x) \mathbb{E}_0\big[\Delta_N(X(t)) \mathbb{1}_{\{X(t) \in W_N\}}\big], \qquad t \to \infty. \quad (4.14)$$

**Remark 4.10.** (i) Relation (4.13) is remarkable since it provides a host of examples of orthogonal polynomial ensembles that appear as the marginal distribution of $h$-transformed diffusions with $h = \Delta_N$ (recall that $\mathbb{P}_0(X(t) \in \mathrm{d}y)$ is a product measure). Explicit examples are the Hermite ensemble for Brownian motion and the Laguerre ensemble for squared Bessel processes, where $\mathbb{P}_0(X(t) \in \mathrm{d}y)$ is the Gamma distribution. Most of the other examples covered by Lemma 4.9 do not seem to be explicit.

(ii) Relation (4.7) may be deduced from (4.14), if the right hand side is asymptotically equivalent when $t$ is replaced by $t - s$ for some $s > 0$. This has not been worked out yet in general, but can be easily seen in a couple of special cases. It would justify the notion 'non-colliding diffusion' for $h$-transformed diffusions with $h = \Delta_N$.

(iii) A natural question is what examples (besides the Hermite ensemble, i.e., Brownian motions; see Section 4.1) lead to processes that can be represented as eigenvalue processes for suitable matrix-valued diffusions. We mention here the Laguerre process, the non-colliding version of squared Bessel processes, which is in distribution equal to the eigenvalue process of a natural processes of complex Wishart matrices ([KO01]; see Remark 2.2(v)). We recall that the real-matrix case, which does not seem to admit an $h$-transform interpretation, is worked out in [Br91].

(iv) Further important examples with physical relevance are derived in [KT04]; in fact, process versions of *all* ten classes of Gaussian random matrices mentioned at the beginning of Section 2 are analysed, and their eigenvalue processes are characterised in terms of non-colliding diffusions.

(v) In [KNT04], independent Brownian motions are conditioned on non-collision up to a *fixed* time, $S$. The result is a time-inhomogeneous diffusion whose transition probabilities depend on $S$. This conditioned process converges towards Dyson's Brownian motions as $S \to \infty$. In [KT03], the distribution of the conditional process is identified in terms of a certain eigenvalue diffusion of a matrix-valued diffusion. Indeed, let $(M_1(t))_{t\geq 0}$ be a Hermitian matrix-valued diffusion whose sub-diagonal and diagonal entries are $\frac{1}{2}N(N+1)$ independent standard real Brownian motions, and let $(M_2(t))_{t\geq 0}$ be an antisymmetric matrix-valued diffusion whose sub-diagonal entries are $\frac{1}{2}N(N-1)$ real independent Brownian bridges (i.e., Brownian motions conditioned on being back to the origin at time $S$). Then the eigenvalue process for the matrix $M_1(t) + \mathrm{i}\, M_2(t)$ is a realisation of the above conditioned Brownian motion process, given that no collision happens by time $S$. The matrix diffusion $(M_1(t) + \mathrm{i}\, M_2(t))_{t\in[0,S]}$ is a one-parameter interpolation between GUE and GOE (hence it is sometimes called a *two-matrix model*). Indeed, recall the well-known independent decomposition of a Brownian motion $(B(t))_{t\geq 0}$ into the Brownian bridge $(B(t) - \frac{t}{S}B(S))_{t\in[0,S]}$ and the linear function $(\frac{t}{S}B(S))_{t\in[0,S]}$ and decompose $M_1(t)$ in that way. Collecting the bridge parts of $M_1(t) + \mathrm{i}\, M_2(t)$ in one process and the remaining variables in the other, we obtain the interpolation.

(vi) *Infinite* systems of non-colliding random processes are considered in [Ba00] and in [KNT04]. The nearest-neighbor discrete-time case is the subject of [Ba00] where the limiting distribution at time $N$ of the left-most walker is derived, conditional on a certain coupling of the total number of left-steps among all the walkers with $N$; the outcome is a certain elementary transformation of the Tracy-Widom distribution for GUE. In [KNT04], a system of $N$ Brownian motions, conditional on non-collision until a fixed time $S$, is analysed in the limit $N \to \infty$ and $S \to \infty$, coupled with each other in various ways. ◇

*4.4.2. The discrete case.*

We present three examples of conditioned random walks on $\mathbb{Z}^N$: the binomial random walk (leading to the Krawtchouk ensemble), the Poisson random walk (leading to the Charlier ensemble) and its de-Poissonized version, the multinomial walk.

For $i = 1, \ldots, N$, let $X_i = (X_i(n))_{n\in\mathbb{N}_0}$ be the *binomial walk*, i.e., at each discrete time unit the walker makes a step of size one with probability $p \in (0,1)$ or stands still otherwise. The walks $X_1, \ldots, X_N$ are assumed independent. Under $\mathbb{P}_x$, the $N$-dimensional process $X = (X_1, \ldots, X_N)$ starts at $X_0 = x \in \mathbb{N}_0^N$. The $\Delta_N$-transformed process on $\mathbb{Z}^N \cap W_N$ has the transition probabilities

$$\widehat{\mathbb{P}}_x(X(n) = y) = \mathbb{P}_x(X(n) = y, T > n)\frac{\Delta_N(y)}{\Delta_N(x)}, \qquad x,y \in \mathbb{Z}^N \cap W_N, n \in \mathbb{N}. \tag{4.15}$$

This marginal distribution, when the process is started at the particular site $x^* = (0,1,2,\ldots,N-1)$, is identified in terms of the Krawtchouk ensemble in (3.8) as follows.

**Lemma 4.11 ($\Delta_N$-transformed binomial walk, [KOR02]).** *Let $x^* = (0,1,2,\ldots,N-1)$. Then, for any $n \in \mathbb{N}$, and $y \in \mathbb{Z}^N \cap W_N$,*

$$\widehat{\mathbb{P}}_{x^*}(X(n) = y) = \mathrm{Kr}_{N,n+N-1,p}(y). \tag{4.16}$$

Such an identification is known only for the particular starting point $x^*$. The proof is based on the Karlin-McGregor formula and some elementary calculations for certain determinants.

The *Poisson random walk*, $X_i = (X_i(t))_{t \geq 0}$, on $\mathbb{N}_0$ makes steps of size one after independent exponential random times. If $X_1, \ldots, X_N$ are independent, the process $X = (X_1, \ldots, X_N)$ on $\mathbb{N}_0$ makes steps after independent exponential times of parameter $N$, and the steps are uniformly distributed on the set of the $N$ unit vectors. The embedded discrete-time walk is the so-called *multinomial walk*; at times $1, 2, 3, \ldots$, a randomly picked component makes a unit step. Lemma 4.4(ii) applies also here, and we may consider the $\Delta_N$-transformed version, both in continuous time and in discrete time. The marginal distribution of the discrete-time process is given in (4.15), and the same formula holds true for the continuous-time version with $n \in \mathbb{N}$ replaced by $t > 0$.

Analogously to the binomial walk, the marginal distributions of both conditioned walks, when the process is started at $x^* = (0, 1, 2, \ldots, N-1)$, may be identified in terms of well-known ensembles, which we introduce first. The *Charlier ensemble* with parameter $\alpha > 0$ and $N \in \mathbb{N}$ is given as

$$\mathrm{Ch}_{N,\alpha}(x) = \frac{1}{Z_{\alpha,N}} \Delta_N(x)^2 \prod_{i=1}^N \frac{\alpha^{x_i}}{x_i!}, \qquad x \in \mathbb{N}_0^N \cap W_N. \qquad (4.17)$$

The *de-Poissonized Charlier ensemble* is defined as

$$\mathrm{dPCh}_{N,n}(x) = \frac{1}{Z_{N,n}} \Delta_N(x)^2 \mathrm{Mu}_{N,n}(x), \qquad x \in \mathbb{N}_0^N \cap W_N, n \in \mathbb{N}_0, , \qquad (4.18)$$

where

$$\mathrm{Mu}_{N,n}(x) = \begin{cases} N^{-n}\binom{n}{x_1,\ldots,x_N} & \text{if } x_1 + \cdots + x_N = n, \\ 0 & \text{otherwise.} \end{cases} \qquad (4.19)$$

Then the free multinomial random walk has the marginals $\mathbb{P}_x(X(n) = y) = \mathrm{Mu}_{N,n}(y - x)$.

**Lemma 4.12 (Conditioned Poisson and multinomial walks, [KOR02]).** *Let $x^* = (0, 1, 2, \ldots, N-1)$.*

(i) *Let $X = (X(t))_{t \geq 0}$ be the Poisson walk, then the marginal distribution of the conditional process satisfies, for any $t > 0$ and $x \in \mathbb{Z}^N \cap W_N$,*

$$\widehat{\mathbb{P}}_{x^*}(X(t) = x) = \mathrm{Ch}_{N,t}(x). \qquad (4.20)$$

(i) *Let $X = (X(n))_{n \in \mathbb{N}_0}$ be the multinomial walk, then the marginal distribution of the conditional process satisfies, for any $n \in \mathbb{N}_0$ and $x \in \mathbb{N}_0^N \cap W_N$,*

$$\widehat{\mathbb{P}}_{x^*}(X(n) = x) = \mathrm{dPCh}_{N,n+N(N-1)/2}(x). \qquad (4.21)$$

The proofs of Lemma 4.12 are based on the Karlin-McGregor formula and explicit calculations for certain determinants.

**Acknowledgement.** The support of the Transregio Sonderforschungsbereich 12 *Symmetries and Universality in Mesoscopic Systems* at Ruhr-Universität

Bochum, where part of this text was written, is gratefully acknowledged. I also thank the Deutsche Forschungsgemeinschaft for awarding a Heisenberg grant, which was realized in 2003/04. Furthermore, I am grateful for the support by, hints from and discussions with some of the experts in the field, which kept me update with the latest developments and helped me to put things in the right respect, among which P. Ferrari, M. Katori, T. Kriecherbauer, N. O'Connell, M. Prähofer, A. Soshnikov, H. Spohn, and C. Tracy.